\documentstyle[amscd]{amsart}
\numberwithin{equation}{section}
\theoremstyle{plain}
\newtheorem{thm}{Theorem}[section]

\newtheorem{lem}[thm]{Lemma}
\newtheorem{prop}[thm]{Proposition}
\begin{document}
\title{Cuntz-Krieger algebras and a generalization of Catalan numbers }
\author{Kengo Matsumoto}
\address{ 
Department of Mathematical Sciences, 
Yokohama City  University,
22-2 Seto, Kanazawa-ku, Yokohama, 236-0027 Japan}
\email{kengo@@yokohama-cu.ac.jp}
\maketitle
\begin{abstract}
We first observe that the relations of the canonical generating isometries 
of the Cuntz algebra ${\cal O}_N$
are naturally related to the $N$-colored Catalan numbers. 
For a directed graph $G$, we generalize the Catalan numbers by using the 
canonical generating partial isometries of the Cuntz-Krieger algebra 
${\cal O}_{A^G}$ for the transition matrix $A^G$ of $G$.
The generalized Catalan numbers $c_n^G, n=0,1,2,\dots$
enumerate the number of Dyck paths and oriented rooted trees for the graph $G$.
Its generating functions will be studied.
\end{abstract}

\def\Zp{{ {\Bbb Z}_+ }}
\def\U{{ {\cal U} }}
\def\S{{ {\cal S} }}
\def\M{{ {\cal M} }}
\def\OL{{ {\cal O}_\Lambda }}
\def\A{{ {\cal A} }}
\def\Ext{{{\operatorname{Ext}}}}
\def\Im{{{\operatorname{Im}}}}
\def\Hom{{{\operatorname{Hom}}}}
\def\Ker{{{\operatorname{Ker}}}}
\def\dim{{{\operatorname{dim}}}}
\def\id{{{\operatorname{id}}}}
\def\OLF{{{\cal O}_{{\frak L}^{Ch(D_F)}}}}
\def\OLN{{{\cal O}_{{\frak L}^{Ch(D_N)}}}}
\def\OLA{{{\cal O}_{{\frak L}^{Ch(D_A)}}}}
\def\LCHDA{{{{\frak L}^{Ch(D_A)}}}}
\def\LCHDN{{{{\frak L}^{Ch(D_N)}}}}
\def\LCHDF{{{{\frak L}^{Ch(D_F)}}}}
\def\LCHLA{{{{\frak L}^{Ch(\Lambda_A)}}}}
\def\LWA{{{{\frak L}^{W(\Lambda_A)}}}}




\section{Introduction}
In the theory of combinatorics,
the process of enumerating objects of various nature  
has been considered (cf. \cite{Ri2}, \cite{Sta}, etc.).
It yields a sequence of positive integers and its generating function.
The Catalan numbers are one of typical examples of such sequences.
The numbers enumerate various objects, 
brackets, Dyck paths, rooted trees, triangulation of polygons, etc.
(cf. \cite{Deu}, \cite{La}, etc. ).    

In this paper, 
we will 
first observe that the operator relations of the canonical generating isometries of the Cuntz algebra ${\cal O}_N$
are naturally related to the $N$-colored Catalan numbers
$c^{(N)}_n,n=0,1,2,\dots.$
The Cuntz-Krieger algebras are natural generalization of the Cuntz algebras
from the view point of topological Markov shifts.
They are defined by a directed graph and have generating partial isometries satisying certain operator relations coming from the structure of the graph.
For a directed graph $G$, we generalize the Catalan numbers by using the 
canonical generating partial isometries of the Cuntz-Krieger algebra 
${\cal O}_{A^G}$ for the transition matrix $A^G$ of $G$.
We call the generalized Catalan numbers ${c_n^G}, n=0,1,2,\dots$
$G$-Catalan numbers.
We will then show that the generalized Catalan numbers enumerate 
the Dyck paths and the oriented rooted trees associated to the graph $G$. 
Let $\{v_1,\dots,v_N\}$ be the vertex set of $G$.
The generalized Catalan numbers ${c_n^G}, {n=0, 1,2,\dots}$
are of the form
$$
c_n^G  = \sum_{i=1}^N c_n^G(i),\qquad n=0,1,2,\dots
$$
where
$c_n^G(i)$ enumerate the numbers rooted at the vertex $v_i$,
and 
$c_0^G(i)$ is defined to be $1$ for $i=1,\dots,N$ so that 
$c_0^G = N$.
They satisfy the  following relation:
\begin{equation}
c_{n+1}^G(i) =\sum_{k=0}^n c_{n-k}^G(i) \sum_{j=1}^N A_G(j,i)c_k^G(j)
\end{equation}
where 
 $A_G(j,i)$ for $v_i, v_j \in V$ denotes the number of directed edges from $v_j$ to $v_i$ in the graph $G$.
Let 
$f^G(x) $ 
be the generating function for the $G$-Catalan numbers
$c_n^G, n=0,1,\dots$.
It is defined by
\begin{equation}
f^G(x) = \sum_{n=0}^{\infty} c_n^G x^n.\label{eqn:fig}
\end{equation}
Let 
$f_i^G(x) $ 
be the generating function for the sequence
$c_n^G(i), n=0,1,\dots$.
The functions satisfy the following equations
\begin{align}
f^G(x) & = \sum_{i=1}^N f_i^G(x), \\
f_i^G(x) & = 1 + x f_i^G(x) \sum_{j=1}^N A_G(j,i) f_j^G(x) 
\qquad \text{ for } i=1,\dots,N.
\end{align}
We will prove that the family 
$f^G_i(x), i=1,\dots,N$
of functions 
is uniquely determined by the above relations by the implicit function theorem
(Theorem 6.5).
The radius of convergence of $f^G_i(x)$ does not depend on $i=1,\dots,N$,
and is determined algebraically as a solution of an eigenvalue problem of a certain matrix 
associated with the transition matrix $A_G$ of the graph $G$ (Theorem 7.2).
We will prove a formula of $c_n^G(i)$ by using the matrix $A_G$. 
Put 
\begin{align*}
F_i(w_1,\dots,w_N) & = (w_i +1) \sum_{j=1}^N A_G(j,i) (w_j +1),\\
F_i^n(w_1,\dots,w_N) & = F_i(w_1,\dots,w_N)^n,\qquad i=1,\dots,N.
\end{align*}
We will show the following integral formulae 
\begin{equation*}
c_n^G(i) =  \frac{1}{2\pi n \sqrt{-1} } \int_C \! \! 
\frac{F_j^n(w_1,\dots,w_N)}{w_j^n} dw_i,\qquad i,j =1,\dots,N
\end{equation*}
hold, where the above integral is a contour integral along a closed curve $C$
around the origin (Theorem 8.1).
By using the above formulae, we will compute some examples.
In particular, the $G$-Catalan numbers for the graph $G$ bellow are computed as
\begin{figure}[htbp]
\begin{center}
\unitlength 0.1in
\begin{picture}( 21.9100,  6.9000)( 16.5000,-24.0600)
%
\special{pn 8}%
\special{ar 2092 2064 78 66  0.0000000 6.2831853}%
%
\special{pn 8}%
\special{ar 3762 2058 80 66  0.0000000 6.2831853}%
%
\special{pn 8}%
\special{ar 1832 2058 182 154  0.3190696 6.0756891}%
%
\special{pn 8}%
\special{pa 2002 2112}%
\special{pa 2006 2106}%
\special{fp}%
\special{sh 1}%
\special{pa 2006 2106}%
\special{pa 1958 2158}%
\special{pa 1982 2154}%
\special{pa 1994 2176}%
\special{pa 2006 2106}%
\special{fp}%
%
\special{pn 8}%
\special{ar 2946 2128 818 278  6.2480092 6.2831853}%
\special{ar 2946 2128 818 278  0.0000000 3.1562795}%
%
\special{pn 8}%
\special{pa 3762 2122}%
\special{pa 3762 2118}%
\special{fp}%
\special{sh 1}%
\special{pa 3762 2118}%
\special{pa 3742 2186}%
\special{pa 3762 2172}%
\special{pa 3782 2186}%
\special{pa 3762 2118}%
\special{fp}%
%
\special{pn 8}%
\special{pa 2134 2006}%
\special{pa 2136 1974}%
\special{pa 2148 1946}%
\special{pa 2166 1918}%
\special{pa 2188 1896}%
\special{pa 2214 1876}%
\special{pa 2240 1858}%
\special{pa 2268 1842}%
\special{pa 2296 1830}%
\special{pa 2326 1816}%
\special{pa 2356 1804}%
\special{pa 2386 1794}%
\special{pa 2416 1784}%
\special{pa 2448 1776}%
\special{pa 2478 1768}%
\special{pa 2510 1760}%
\special{pa 2542 1754}%
\special{pa 2572 1748}%
\special{pa 2604 1742}%
\special{pa 2636 1738}%
\special{pa 2668 1734}%
\special{pa 2700 1730}%
\special{pa 2732 1728}%
\special{pa 2762 1724}%
\special{pa 2794 1722}%
\special{pa 2826 1720}%
\special{pa 2858 1720}%
\special{pa 2890 1718}%
\special{pa 2922 1716}%
\special{pa 2954 1716}%
\special{pa 2986 1718}%
\special{pa 3018 1718}%
\special{pa 3050 1718}%
\special{pa 3082 1720}%
\special{pa 3114 1722}%
\special{pa 3146 1726}%
\special{pa 3178 1728}%
\special{pa 3210 1730}%
\special{pa 3242 1736}%
\special{pa 3274 1740}%
\special{pa 3306 1744}%
\special{pa 3336 1750}%
\special{pa 3368 1756}%
\special{pa 3400 1764}%
\special{pa 3430 1770}%
\special{pa 3462 1778}%
\special{pa 3492 1788}%
\special{pa 3522 1796}%
\special{pa 3554 1806}%
\special{pa 3582 1820}%
\special{pa 3612 1832}%
\special{pa 3640 1846}%
\special{pa 3668 1862}%
\special{pa 3694 1880}%
\special{pa 3718 1902}%
\special{pa 3740 1926}%
\special{pa 3758 1952}%
\special{pa 3768 1982}%
\special{pa 3768 2000}%
\special{sp}%
%
\special{pn 8}%
\special{pa 2134 1992}%
\special{pa 2134 2006}%
\special{fp}%
\special{sh 1}%
\special{pa 2134 2006}%
\special{pa 2154 1940}%
\special{pa 2134 1954}%
\special{pa 2114 1940}%
\special{pa 2134 2006}%
\special{fp}%
\end{picture}%
\end{center}
\caption{}
\end{figure}
$$
c_n^G = 
\frac{2}{n} \binom{3n+1}{n-1} + \frac{1}{n} \binom{3n}{n-1}
=
\frac{2}{n+1} \binom{3n}{n}.
$$
Hence 
$$
c_0^{G}=2,\quad
c_1^{G}=3,\quad
c_2^{G}=10,\quad
c_3^{G}=42,\quad
c_4^{G}=198,\quad \dots.
$$
This sequence is regarded to be the Fibonacci version of the Catalan numbers.
The radius of convergence of the function 
$f^G(x)$ is computed to be $\frac{4}{27}$.

We will also define generalized Catalan numbers 
$c^{G,\varphi}_n$ by using KMS state $\varphi$ 
on the $C^*$-algebra ${\cal O}_{A^G}$ for the gauge action.
The numbers $c^{G,\varphi}_n$ are computed by 
the $G$-Catalan numbers with 
the Perron-Frobenius eigenvalue and its eigenvector for the matrix $A^G$ (Proposition 9.3).

\section{The Cuntz algebras and Catalan numbers}
Throughout this section $N$ is a fixed positive integer greater than $1$.
For a set $\Omega$, we denote by $|\Omega|$ the number of $\Omega$.  
Let $S_1,\dots, S_N$ be a family of bounded linear operators on a Hilbert space satisfying the following condition
\begin{equation}
\sum_{j=1}^N S_jS_j^* = 1, \qquad S_i^* S_i = 1 
\quad \text{ for } \quad i=1,\dots,N \label{eqn:on}
\end{equation}
Let ${\cal O}_N$ be the $C^*$-algebra generated by the family  
$S_1,\dots, S_N$.
The algebra ${\cal O}_N$ is called the Cuntz algebra (\cite{C}).
It is well-known that the algebraic structure of ${\cal O}_N$ does not depend on the choice of the family $S_1,\dots, S_N$
satisfying the relations  \eqref{eqn:on}.
For a word $\mu = (\mu_1,\dots, \mu_n)$ of $\{ 1,\dots,N\}$,
we put
$$
S_{\mu} = S_{\mu_1}\cdots S_{\mu_n}
\quad
\text{ and }
\quad
S_{\mu}^* = S_{\mu_n}^*\cdots S_{\mu_1}^*.
$$
Put
$\Sigma_N = \{ S_1^*,\dots,S_N^*,S_1,\dots,S_N\}.$
Let
$W_{2n}^N$ be the set of all words of $\Sigma_N$ of length $2n$:
$$
W_{2n}^N = \{ (X_1,\dots,X_{2n}) \mid X_i \in \Sigma_N, \, i=1,\dots,2n \}.
$$
For $X=(X_1,\dots,X_{2n})\in W_{2n}^N$,
define 
$\pi_N (X) $ to be the element 
$X_1\cdots X_{2n}$ in ${\cal O}_N$.
By \cite[1.3 Lemma]{C}, 
every $\pi_N(X)$ for a  word $X$ in $W_{2n}^N$ is one of  the forms:
$$
0, \qquad 1, \qquad  S_{\mu} S_{\nu}^* \quad\text{ for some words }
\mu, \nu \text{ of } \{1,\dots,N\}.
$$ 
We set
$$
B_n^N = \{ X \in W_{2n}^N \mid \pi_N(X) = 1 \}.
$$
Define the numbers
$$
c_0^{(N)}  = N,\qquad
c_n^{(N)} = | B_n^N |\quad  \text{ for } n=1,2,\dots.
$$
The following proposition is a key in our further discussions
\begin{lem}
Let $X = (X_1, \dots, X_{2n})$ be a word in $B_n^N$.
\hspace{3cm}
\begin{enumerate}
\renewcommand{\labelenumi}{(\roman{enumi})}
\item
The number of $S_i^*, i=1,\dots,N$ and $S_i, i=1,\dots,N$
in $X$ is the same. 
\item
The number of $S_i^*, i=1,\dots,N$ in any starting segment from the leftmost of $X$
 is not less than the number of $S_i, i=1,\dots,N$ in the same segment.
\end{enumerate}
\end{lem}
\begin{pf}
(i) is clear.

(ii) Let $\mu,\nu$ be the word of $\{1,\dots,N\}$ with same length.
Then \eqref{eqn:on} implies
$$
S_{\mu}^* S_{\nu} = 
\begin{cases}
1 & \text{ if } \mu = \nu,\\
0 & \text{ if } \mu \ne \nu.
\end{cases}
$$
Hence, if the number of $S_i^*, i=1,\dots,N$ in some starting segment from leftmost of $X$
 is  less than the number of $S_i, i=1,\dots,N$ in the same segment.
The operator $\pi_N(X)$ is reduced to the word of the form 
$S_j \pi_N(Y)$ for some $j=1,\dots,N$ and some word $Y$ in $\Sigma_N$.
The word 
$S_j \pi_N(Y)$ however is not abe  to be $1$ for any word $Y$.
Thus the assertion holds.
\end{pf}
Let $c_n$ be the usual Catalan number $\frac{1}{n+1}\binom{2n}{n}$.
The above two properties in Lemma 2.1 characterize the regular bracket structures (\cite[p.26]{La}),
so that we have
\begin{prop}
$$
c_n^{(N)} = N^n \times c_n, \qquad n=0,1,\dots.
$$
The numbers 
$ N^n \times c_n$ are called the $N$-colored Catalan numbers.
\end{prop}
\begin{pf}
Let 
$(_1,\dots,(_N, )_1,\dots,)_N$
be the $N$ pairs of $2N$ brackets.
Consider the following correspondence between
the brackets $(_1,\dots,(_N, )_1,\dots,)_N$ 
and the operators
$S_1^*,\dots,S_N^*,S_1,\dots,S_N$
such as
$$(_i \rightarrow S_i^*,\qquad
)_i \rightarrow S_i \qquad \text{ for } i=1,\dots,N.
$$
One easily sees that the set of regular brackets exactly corresponds to the set of words of $B_n^N$.
\end{pf}
 For example, the word
$$
S_{f_4}^*S_{f_3}^*S_{f_3} S_{f_2}^*S_{f_1}^*
S_{f_1}S_{f_2}S_{f_4} \cdot S_{f_3}^* S_{f_2}^*S_{f_2}S_{f_3}
$$
corresponds to the word of brackets
$$
(_{f_4}  (_{f_3}  )_{f_3} (_{f_2}  (_{f_1}  
)_{f_1} )_{f_2} )_{f_4} \cdot (_{f_3} (_{f_2} )_{f_2} )_{f_3}
$$
for $f_1,f_2,f_3,f_4 \in \{1,\dots,N\}$.
We note that for every word $X = X_1,\dots X_{2n}$ of $B_n^N$,
the leftmost symbol $X_1$ must be $S_i^*$ for some $i$
and the rightmost symbol $X_{2n}$ must be $S_j$ for some $j$.
We denote  $i$ by $l(X)$ and $j$ by $r(X)$ respectively.

\section{The Cuntz-Krieger algebras and generalized Catalan numbers}

Let $G =(V,E)$ be a directed finite graph with vertex set $V$
and edge set $E$.
Consider the associated edge matrix $A^G$ for $G$ defined by for 
$e,f \in E$
$$
A^G(e,f)
=\begin{cases}
1 & \text{ if } t(e) = s(f),\\
0 & \text{ otherwise}
\end{cases}
$$
where $t(e), s(f)$ denote the terminal vertex 
of edge $e$ and the source vertex of edge $f$.
We henceforth assume that every vertex of $G$ has an incoming edge and an outgoing edge so that $A^G$ has no zero rows or columns.
Consider the Cuntz-Krieger algebra ${\cal O}_{A^G}$ for the matrix $A^G$
that is the universal $C^*$-algebra generated by 
partial isometries $S_e, e \in E$ subject to the following relations:
\begin{equation}
\sum_{f \in E} S_f S_f^* = 1, 
\qquad
S_e^* S_e = \sum_{f \in E} A^G(e,f) S_fS_f^* \quad \text{ for } e \in E.
\label{eqn:oa}
\end{equation}
For a vertex $v \in V$ we put
$$
P_v = \sum_{f \in E, v=s(f)} S_f S_f ^*.
$$
By the above relations, we know that
\begin{enumerate}
\renewcommand{\labelenumi}{(\roman{enumi})}
\item
$P_v P_u =0$ if $u\ne v$ for $u,v \in V$,
\item
$ P_{t(e)} = S_e^* S_e $ for $e \in E$.
\end{enumerate}
The algebra ${\cal O}_{A^G}$ is often called a graph algebra 
and denoted by ${\cal O}_G$ \cite{KPRR,KPR}.

Let $\{e_1, \dots, e_{|E|}\}$ be the edge set $E$. 
Put  
$\Sigma_G = \{ S_{e_1}^*,\dots,S_{e_{|E|}}^*,S_{e_1},\dots,S_{e_{|E|}} \}$.
We denote by $W_{2n}^G$ the set of all words of 
$\Sigma_G$ of length $2n$. 
For a word $X = (X_1,\dots,X_{2n}) \in W_{2n}^G,$
define
$\pi_G(X) = X_1 \cdots X_{2n} \in {\cal O}_{A^G}$
as an element of the algebra ${\cal O}_{A^G}$.
Let $S_1,\dots, S_{|E|}$ 
be the canonical generating isometries of 
the Cuntz-algebra ${\cal O}_{|E|}$ satisfying \eqref{eqn:on}.
Put 
$
\Sigma_{|E|} = \{ S_1^*,\dots,S_{|E|}^*,S_1,\dots,S_{|E|} \}
$
and define the correspondence 
$$
\Phi_G :\Sigma_G
 \longrightarrow \Sigma_{|E|}
$$
by setting
$$ 
\Phi_G(S_{e_i}^*) = S_i^*,\qquad 
\Phi_G(S_{e_i}) = S_i,  \qquad i=1,\dots,{|E|}.
$$
The correspondence $\Phi_G$ is naturally extended to words of 
$\Sigma_G$.
We set
$$
B_n^G = \{ X \in W_{2n}^G 
\mid \pi_G(X) \ne 0, \Phi_G(X) \in B_n^{|E|} \}.
$$
A word $X =(X_1,\dots, X_{2n})$ of $B_n^G$
is called a $G$-{\it Catalan word}.
Hence a word  
$X =(X_1,\dots, X_{2n})$ of $\Sigma_G$
is a $G$-Catalan word if and only if
$X_1\cdots X_{2n} \ne 0$ in ${\cal O}_{A^G}$
and  
$\Phi_G(X_1)\cdots \Phi_G(X_{2n}) =1 $ in ${\cal O}_{|E|}$.
We set 
$$
c_0^{G} =N, \qquad 
c_n^{G} = | B_n^G | \quad \text{ for } n= 1,2, \dots
$$
where $N$ denotes the number $| V |$ of the vertex set $V$.
We call the sequaence 
$c_n^{G}, n=0,1,\dots $ 
{\it the generalized Catalan number associated with the graph } $G$,
or $G$-{\it Catalan numbers } for brevity.

Let $\Lambda_G$ be the topological Markov shift 
$$
\Lambda_G = \{ {(f_i)}_{i \in \Bbb Z}\in E^{\Bbb Z}
 \mid A^G(f_i,f_{i+1}) = 1, i \in \Bbb Z \}
$$ 
defined by the matrix $A^G$.
We denote by $\Lambda_G^*$ 
the set of all admissible words of the subshift $\Lambda_G$ (cf. \cite{LM}).
\begin{lem}
For $ f_1,\dots,f_k \in E$, 
the following identity
$$
S_{f_k}^* S_{f_{k-1}}^* \cdots S_{f_2}^*
S_{f_1}^* S_{f_1} S_{f_2} \cdots S_{f_{k-1}} S_{f_k}
=
A^G(f_1,f_2)A^G(f_2,f_3)\cdots A^G(f_{k-1},f_k ) S_{f_k}^* S_{f_k}
$$
holds.
If in particular $f_1 \cdots f_k \in \Lambda_G^*$, we have
$$
S_{f_k}^* S_{f_{k-1}}^* \cdots S_{f_2}^* S_{f_1}^* 
S_{f_1} S_{f_2} \cdots S_{f_{k-1}} S_{f_k}
=S_{f_k}^* S_{f_k}.
$$
\end{lem}
\begin{pf}
By using relations \eqref{eqn:oa} recursively,
the above identities are straightforward.
\end{pf}

\begin{lem}
For every $X$ in $B_n^G$,
there exists a vertex $v(X) \in V$ such that
$\pi_G(X) = P_{v(X)}.$
\end{lem}
\begin{pf}
By  Proposition 2.2 and Lemma 3.1, one sees that 
for every $X$ in $B_n^G$, the element $\pi_G(X)$ is of the form
$$
\pi_G(X) = 
S_{f_{k_1}}^* S_{f_{k_1}}
S_{f_{k_2}}^*S_{f_{k_2}}
\cdots 
S_{f_{k_m}}^*S_{f_{k_m}}
$$
for some words
$
f_{k_1}, f_{k_2},
\cdots f_{k_m}.
$
Since 
$S_{f_{k_j}}^*S_{f_{k_j}} = P_{t(f_{k_j})}$,
one has
$$
\pi_G(X) = P_{t(f_{k_1})} P_{t(f_{k_2})} \cdots P_{t(f_{k_m})}.
$$
As $\pi_G(X) \ne 0$ and $P_u P_v = 0 $ for $u \ne v$,
one obtains 
$$
t(f_{k_1})= t(f_{k_2}) = \cdots = t(f_{k_m}).
$$
By putting
$v = t(f_{k_1})$,
one concludes that
$\pi_G(X) = P_v$.
\end{pf}

Put for $e \in E$
$$
B_n^G[e]
= \{ X=(X_1,\dots,X_{2n}) \in B_n^G \mid 
(S_e^*,X_1,\dots,X_{2n}, S_e ) \in B_{n+1}^G \}.
$$
We note the following lemma 
\begin{lem}
A word $X$ in $B_n^G$ belongs to $
B_n^G[e]$ if and only if $v(X) = s(e)$.
\end{lem}
\begin{pf}
For $ X=(X_1,\dots,X_{2n}) \in B_n^G $,
it follows that
$$
\pi_G((S_e^*,X_1,\dots,X_{2n}, S_e )) = S_e^* \pi_G(X) S_e = S_e^* P_{v(X)} S_e.
$$
As 
$S_e^* P_{v(X)} S_e \ne 0$ if and only if 
$P_{v(X)} \ge S_e S_e^*$.
The latter condition is equivalent to the condition
$v(X) = s(e)$.
Hence
the word
$(S_e^*,X_1,\dots,X_{2n}, S_e )$ belongs to
$ B_{n+1}^G$ if and only if
$v(X) = s(e)$.
\end{pf}

For
 $ X=(X_1,\dots,X_{2n}) \in B_n^G $,
 the leftmost symbol $X_1$ must be $S_e^*$ for some $e \in E$ 
 and the rightmost symbol $X_{2n}$ must be $S_f$ for some $f \in E$.
We denote  $e$ by $l(X)$ and $f$ by $r(X)$ respectively.
Hence $X$ is of the form
$$
X = (S_{l(X)}^*, X_2, \dots X_{2n-1}, S_{r(X)}).
$$
As $S_{l(X)}^*S_{l(X)} =S_{r(X)}^*S_{r(X)} = P_{v(X)},$
one sees $t(l(X)) = t(r(X)) =v(X).$
The following property for brackets is well-known.
\begin{lem}
For a $G$-Catalan word 
 $X =(X_1,\dots, X_{2n+2})$ in $B_{n+1}^G$,
there uniquely exists 
$k \in \Bbb N$ with $0 \le k \le n$
such that
$$
(X_2,\dots,X_{2k+1}) \in B_k^G, 
\quad
X_{2k+2} = X_1^*,\quad
(X_{2k+3},\dots,X_{2n+2}) \in B_{n-k}^G.
$$
\end{lem}
This lemma means by putting
$$
Y=(X_2,\dots,X_{2k+1}), \qquad
Z =(X_{2k+3},\dots,X_{2n+2}),
$$
the word $X$ is decomposed as 
$$
X = (S_{l(X)}^*, Y, S_{l(X)}, Z ), \qquad
Y \in B_k^G, \qquad
Z \in B_{n-k}^G,
$$
in a unique way.
\begin{lem}
For 
$e,f \in E$,
and
$Y =(Y_1,\dots,Y_{2k}) \in B_k^G, 
Z =(Z_1,\dots, Z_{2(n-k)}) \in B_{n-k}^G,
$
the word
$(S_f^*, Y_1,\dots,Y_{2k}, S_f, Z_1,\dots, Z_{2(n-k)}) \in W_{2n+2}^G
$
belongs to
$B_{n+1}^G[e]
$
if and only if
$$
Y \in B_k^G[f], \quad
Z  \in B_{n-k}^G[e] \quad
\text{ and }\quad
t(f) = s(e).
$$
\end{lem}
\begin{pf}
By Lemma 3.3, 
$ Y$ belongs to $ B_k^G [f] $ 
if and only if
$v(Y) = s(f),$
and 
$Z$ belongs to $ B_{n-k}^G [e]$ if and only if 
$ v(Z)= s(e) $.
As we have
$$
\pi_G((S_f^*, Y_1,\dots,Y_{2k}, S_f, Z_1,\dots, Z_{2(n-k)}))
= S_f^* P_{v(Y)} S_f P_{v(Z)},
$$
the above element is not zero if and only if 
$v(Y) = s(f)$ and $v(Z) =t(f)$.
The latter condition is equivalent to the conditions
$$
Y \in B_k^G[f], \quad
Z  \in B_{n-k}^G[e] \quad
\text{ and }\quad
t(f) = s(e).
$$
\end{pf}
Hence we have 
\begin{lem} \label{lem:klem}
For $e \in E$, the equality 
$$
B_{n+1}^{G}[e] = 
 \bigsqcup\begin{Sb}
 f \in E \\ A^G(f,e) =1
\end{Sb}
\bigsqcup^n_{k=0} \, B_{k}^G [f] \times B_{n-k}^G [e]
$$
holds through the correspondence
$$
X = (S_{f}^*, Y, S_f, Z )\in B_{n+1}^G[e]
\longrightarrow 
(Y, Z) \in B_{k}^{G}[f] \times B_{n-k}^{G}[e].
$$ 
\end{lem}
We set
$$
c_0^G[e] = 1, \qquad
c_n^G[e] = |B_n^G[e]| \quad \text{ for }\quad n=1,2,\dots, \quad e \in E.
$$
For $e,f \in E$, 
define an equivalence relation 
$e \sim f $ by the condition $s(e) = s(f)$.
The equivalence relation
$e\sim f$ implies 
$B_n^G[e]= B_n^G[f]$ 
and hence $c_n^G[e] = c_n^G[f]$.
For a vertex $u \in V$, we define
$$
B_n^G(u) = \{ X \in B_n^G \mid v(X) =u \}
$$
so that 
$B_n^G = 
\sqcup_{u \in V}
B_n^G(u).
$
We set for a vertex $u \in V$,
$$
c_0^G(u) = 1, \qquad
c_n^G(u) = |B_n^G(u)| \quad 
\text{ for } n=1,2,\dots, \quad u \in V.
$$
Therefore we have 
\begin{prop}
\hspace{6cm}
\begin{enumerate}
\renewcommand{\labelenumi}{(\roman{enumi})}
\item
$c_n^G = \sum\limits_{u \in V} c_n^G(u).$
\item
$
c_{n+1}^G[e] =\sum\limits_{k=0}^n c_{n-k}^G[e] 
\sum\limits_{f \in E} 
A^G(f,e)c_k^G[f].
$
\item
 if $s(e) = u,$ we have
$
c_n^G(u) = c_n^G[e].
$
\end{enumerate}
\end{prop}
\begin{pf}
The assertions are all obvious.
\end{pf}

\section{Dyck paths  associated with the graph $G$}
We will define Dyck paths associated with a given directed graph $G$.
We will enumerate them and define the sequence
$
d_n^G, n=0,1,\dots
$
of numbers.
We will prove that
$$ d_n^G = c_n^G \qquad \text{ for } n=0,1,\dots.
$$ 
Let $G=(V,E)$ be a directed graph.
Let $G^*$ denote the transposed graph of $G$.
The vertex set $V^*$ of $G^*$ is $V$ 
and the edge set $E^*$ of $G^*$ is the edges reversing the directions of edges of $G$.

A Dyck path $\gamma$ is a continuous broken line 
located in the upper half plane and 
consisting of vectors $(1,1)$ and $(1,-1)$  
starting at the origin and ending at the $x$-axis (see Figure 2).
\begin{figure}[htbp]
\begin{center}
\unitlength 0.1in
\begin{picture}( 60.4000, 14.5000)(  3.7000,-17.0000)
%
\special{pn 20}%
\special{pa 370 1660}%
\special{pa 930 1100}%
\special{fp}%
\special{pa 930 1100}%
\special{pa 1210 1380}%
\special{fp}%
\special{pa 1210 1380}%
\special{pa 1770 820}%
\special{fp}%
\special{pa 1770 820}%
\special{pa 2610 1660}%
\special{fp}%
\special{pa 2610 1660}%
\special{pa 3730 540}%
\special{fp}%
\special{pa 3730 540}%
\special{pa 4290 1100}%
\special{fp}%
\special{pa 4290 1100}%
\special{pa 4570 820}%
\special{fp}%
\special{pa 4570 820}%
\special{pa 5410 1660}%
\special{fp}%
%
\special{pn 20}%
\special{pa 370 1660}%
\special{pa 5380 1670}%
\special{fp}%
%
\special{pn 4}%
\special{pa 650 280}%
\special{pa 650 1660}%
\special{dt 0.027}%
%
\special{pn 4}%
\special{pa 930 260}%
\special{pa 940 1650}%
\special{dt 0.027}%
%
\special{pn 4}%
\special{pa 1210 260}%
\special{pa 1210 1640}%
\special{dt 0.027}%
%
\special{pn 4}%
\special{pa 1490 260}%
\special{pa 1490 1690}%
\special{dt 0.027}%
%
\special{pn 4}%
\special{pa 1770 250}%
\special{pa 1770 1680}%
\special{dt 0.027}%
%
\special{pn 4}%
\special{pa 2050 270}%
\special{pa 2050 1700}%
\special{dt 0.027}%
%
\special{pn 4}%
\special{pa 2330 270}%
\special{pa 2330 1700}%
\special{dt 0.027}%
%
\special{pn 4}%
\special{pa 2610 260}%
\special{pa 2610 1690}%
\special{dt 0.027}%
%
\special{pn 4}%
\special{pa 2890 260}%
\special{pa 2890 1690}%
\special{dt 0.027}%
%
\special{pn 4}%
\special{pa 3170 260}%
\special{pa 3170 1690}%
\special{dt 0.027}%
%
\special{pn 4}%
\special{pa 3450 260}%
\special{pa 3450 1690}%
\special{dt 0.027}%
%
\special{pn 4}%
\special{pa 3730 260}%
\special{pa 3730 1690}%
\special{dt 0.027}%
%
\special{pn 4}%
\special{pa 4010 260}%
\special{pa 4010 1690}%
\special{dt 0.027}%
%
\special{pn 4}%
\special{pa 4290 260}%
\special{pa 4290 1690}%
\special{dt 0.027}%
%
\special{pn 4}%
\special{pa 4570 260}%
\special{pa 4570 1690}%
\special{dt 0.027}%
%
\special{pn 4}%
\special{pa 4850 260}%
\special{pa 4850 1690}%
\special{dt 0.027}%
%
\special{pn 4}%
\special{pa 5130 270}%
\special{pa 5130 1700}%
\special{dt 0.027}%
%
\special{pn 4}%
\special{pa 5410 260}%
\special{pa 5410 1690}%
\special{dt 0.027}%
%
\special{pn 20}%
\special{pa 5680 270}%
\special{pa 5450 290}%
\special{ip}%
\special{pa 6360 290}%
\special{pa 6410 310}%
\special{ip}%
%
\special{pn 4}%
\special{pa 370 260}%
\special{pa 5410 270}%
\special{dt 0.027}%
%
\special{pn 4}%
\special{pa 380 540}%
\special{pa 5420 550}%
\special{dt 0.027}%
%
\special{pn 4}%
\special{pa 370 820}%
\special{pa 5410 830}%
\special{dt 0.027}%
%
\special{pn 4}%
\special{pa 380 1110}%
\special{pa 5420 1120}%
\special{dt 0.027}%
%
\special{pn 4}%
\special{pa 380 1390}%
\special{pa 5420 1400}%
\special{dt 0.027}%
%
\special{pn 4}%
\special{pa 370 270}%
\special{pa 370 1650}%
\special{dt 0.027}%
\end{picture}%
\end{center}
\caption{}
\end{figure}
For a Dyck path 
$\gamma = (\gamma_1,\dots,\gamma_{2n})$,
where
$\gamma_i$ is one of vectors $(1,1)$ and $(1,-1)$,
if $\gamma_i$ is a vector $(1,1)$,
there uniquely exists $\gamma_{i+k}$ satisfying the following conditions:
\begin{enumerate}
\item $\gamma_{i+k}$ is a vector $(1,-1)$.
\item $(\gamma_{i+1}, \gamma_{i+2},\dots,\gamma_{i+k-1})$ is a Dyck path
of length $k-1$ ( hence $k-1$ is even).
\end{enumerate}
 
 We call the edge $\gamma_{i+k}$ the partner of $\gamma_i$.

For an edge $e \in E$, we denote by $e^*$ the edge of $G^*$ 
obtained by reversing the direction of $e$. 
A $G$-{\it Dyck path}\, of length $2n$ is a Dyck path $\gamma$ labeled
$\{ e^*, e \mid e \in E \}$ by the following rules:
\begin{enumerate}
\item vectors $(1,1)$ are labeled $e^*$ for $ e \in E$,
\item vectors $(1,-1)$ are labeled $e$ for $ e \in E$,
\item a vector $(1,1)$ labeled $e^*$ follows a vector $(1,1)$ labeled $f^*$
if and only if $t(f^*) = s(e^*)$ in $G^*$,
\item a vector $(1,1)$ labeled $e^*$ follows a vector $(1,-1)$ labeled $f$
if and only if $t(f) = s(e^*)$ in $G^*$,
\item a vector $(1,-1)$ labeled $e$ follows a vector $(1,1)$ labeled $f^*$
if and only if $e = f$,
\item the partner of a vector $(1,1)$ labeled $e^*$ is labeled by $e$.
\end{enumerate}
Let $D_n^G$ be the set of all $G$-Dyck paths of length $2n$.
For $\gamma = (\gamma_1,\dots,\gamma_{2n}) \in D_n^G$,
the leftmost vector $\gamma_1$ must be $e^*$ for some $e^* \in E^*$,
and  
the rightmost vector $\gamma_{2n}$ must be $f$ for some $f \in E$.
We denote $e^*$ and $ f$ by
$l(\gamma)^*$ and $r(\gamma)$ respectively.
We set
$$
d_0^G  = N, \qquad
d_n^G = |D_n^G| \quad \text{ for } \quad n=1,2,\dots.
$$
The following lemmas are paralle to lemmas in the previous section.
\begin{lem}
For a $G$-Dyck path 
 $\gamma =(\gamma_1,\dots, \gamma_{2n+2})$ in $D_{n+1}^G$,
there uniquely exists 
$k \in \Bbb N$ with $0 \le k \le n$
such that
$$
(\gamma_2,\dots,\gamma_{2k+1}) \in D_k^G, 
\quad
\gamma_{2k+2} = \gamma_1^*,\quad
(\gamma_{2k+3},\dots,\gamma_{2n+2}) \in D_{n-k}^G.
$$
\end{lem}
Put for $e \in E$
$$
D_n^G [e]
= \{ (\gamma_1,\dots,\gamma_{2n}) \in D_n^G \mid 
(e^*,\gamma_1,\dots,\gamma_{2n},e ) \in D_{n+1}^G \}.
$$
\begin{lem}
For 
$e,f \in E$,
and
$\eta =(\eta_1,\dots,\eta_{2k}) \in D_k^G, 
\zeta =(\zeta_1,\dots, \zeta_{2(n-k)}) \in D_{n-k}^G,
$
the path
$(f^*, \eta_1,\dots,\eta_{2k}, f, \zeta_1,\dots, \zeta_{2(n-k)}) $
belongs to
$D_{n+1}^G[e]
$
if and only if
$$
\eta \in D_k^G[f], \quad
\zeta  \in D_{n-k}^G[e] \quad
\text{ and }\quad
t(f) = s(e).
$$
\end{lem}
We set
$$
d_0^G[e] = 1, \qquad
d_n^G[e] = |D_n^G[e]| \quad \text{ for } \quad n=1,2,\dots, \quad e \in E.
$$
Recall that the equivalence relation $\sim$ in $E$ is defined by
$e \sim f $ if $s(e) = s(f)$.
Hence we have
$e\sim f$ implies 
$D_n^G[e]= D_n^G[f]$ and hence $d_n^G[e] = d_n^G[f]$.
For a vertex $u \in V$, we define
$$
D_n^G(u) = \{ \gamma \in D_n^G \mid t(r(\gamma)) =u \}
$$
and
$$
d_0^G(u) = 1, \qquad
d_n^G(u) = |D_n^G[u]| \quad \text{ for } \quad n=1,2,\dots, \quad u \in V.
$$
Therefore we have 
\begin{prop}
\hspace{6cm}
\begin{enumerate}
\renewcommand{\labelenumi}{(\roman{enumi})}
\item
$d_n^G = \sum\limits_{u \in V} d_n^G(u),$
\item
$
d_{n+1}^G[e] =\sum\limits_{k=0}^n d_{n-k}^G[e] 
\sum\limits_{f \in E}
A^G(f,e)d_k^G[f],
$
\item
 if $s(e) = u,$ we have
$
d_n^G(u) = d_n^G[e].
$
\end{enumerate}
\end{prop}
\begin{pf}
The assertions are all obvious.
\end{pf}
Therefore we have
\begin{thm}
For $n=0,1,\dots$,
we have  
$d_n^G[e] = c_n^G[e]$ for $e \in E$
and 
$d_n^G(u)= c_n^G(u)$ for $u \in V$,
so that 
$$
d_n^G = c_n^G. 
$$
\end{thm}

\section{Trees associated with graph $G$}
We will define trees associated with a given directed graph $G$.
Let
$
t_n^G
$
be the numbers of such trees with $n$ edges.
We will prove that
$$ t_n^G = c_n^G \qquad \text{ for } n=0,1,\dots.
$$ 
A rooted tree is a plane tree with a distingushed vertex.
The distinguished vertex is called the root (see Figure 3).
\begin{figure}[htbp]
\begin{center}
\unitlength 0.1in
\begin{picture}( 27.5500, 17.8300)(  2.6200,-22.5000)
%
\special{pn 8}%
\special{sh 0.600}%
\special{ar 1938 640 10 10  0.0000000 6.2831853}%
%
\special{pn 8}%
\special{pa 1938 480}%
\special{pa 1938 480}%
\special{fp}%
%
\special{pn 8}%
\special{sh 0.600}%
\special{ar 1438 1440 10 10  0.0000000 6.2831853}%
%
\special{pn 8}%
\special{sh 0.600}%
\special{ar 2438 1440 10 10  0.0000000 6.2831853}%
%
\special{pn 8}%
\special{sh 0.600}%
\special{ar 1938 2240 10 10  0.0000000 6.2831853}%
%
\special{pn 8}%
\special{sh 0.600}%
\special{ar 1438 2240 10 10  0.0000000 6.2831853}%
%
\special{pn 8}%
\special{sh 0.600}%
\special{ar 938 2240 10 10  0.0000000 6.2831853}%
%
\special{pn 8}%
\special{sh 0.600}%
\special{ar 2438 2240 10 10  0.0000000 6.2831853}%
%
\special{pn 8}%
\special{sh 0.600}%
\special{ar 2938 2240 10 10  0.0000000 6.2831853}%
%
\special{pn 8}%
\special{pa 1938 630}%
\special{pa 1938 480}%
\special{fp}%
%
\special{pn 8}%
\special{ar 1438 1270 10 10  0.0000000 6.2831853}%
%
\special{pn 8}%
\special{pa 1438 1430}%
\special{pa 1438 1280}%
\special{fp}%
%
\special{pn 8}%
\special{pa 1438 1260}%
\special{pa 1438 1110}%
\special{fp}%
%
\special{pn 8}%
\special{pa 938 2230}%
\special{pa 938 2080}%
\special{fp}%
%
\special{pn 8}%
\special{ar 938 2070 10 10  0.0000000 6.2831853}%
%
\special{pn 8}%
\special{pa 938 2060}%
\special{pa 938 1910}%
\special{fp}%
%
\special{pn 8}%
\special{ar 938 1910 10 10  0.0000000 6.2831853}%
%
\special{pn 8}%
\special{ar 1438 1110 10 10  0.0000000 6.2831853}%
%
\special{pn 8}%
\special{pa 1438 2230}%
\special{pa 1438 2080}%
\special{fp}%
%
\special{pn 8}%
\special{pa 1938 2230}%
\special{pa 1938 2080}%
\special{fp}%
%
\special{pn 8}%
\special{pa 2438 2230}%
\special{pa 2438 2080}%
\special{fp}%
%
\special{pn 8}%
\special{pa 2938 2230}%
\special{pa 2938 2080}%
\special{fp}%
%
\special{pn 8}%
\special{pa 938 1900}%
\special{pa 938 1750}%
\special{fp}%
%
\special{pn 8}%
\special{ar 1438 2070 10 10  0.0000000 6.2831853}%
%
\special{pn 8}%
\special{pa 1438 2060}%
\special{pa 1438 1910}%
\special{fp}%
%
\special{pn 8}%
\special{ar 1938 2070 10 10  0.0000000 6.2831853}%
%
\special{pn 8}%
\special{ar 2438 2070 10 10  0.0000000 6.2831853}%
%
\special{pn 8}%
\special{ar 2938 2070 10 10  0.0000000 6.2831853}%
%
\special{pn 8}%
\special{ar 1438 1910 10 10  0.0000000 6.2831853}%
%
\special{pn 8}%
\special{ar 944 1740 10 10  0.0000000 6.2831853}%
%
\special{pn 8}%
\special{pa 1438 2060}%
\special{pa 1628 1920}%
\special{fp}%
%
\special{pn 8}%
\special{pa 1938 2060}%
\special{pa 1938 1910}%
\special{fp}%
%
\special{pn 8}%
\special{pa 2438 2230}%
\special{pa 2518 2080}%
\special{fp}%
\special{pa 2938 2230}%
\special{pa 3008 2086}%
\special{fp}%
\special{pa 2938 2230}%
\special{pa 2878 2090}%
\special{fp}%
%
\special{pn 8}%
\special{pa 2438 1430}%
\special{pa 2518 1280}%
\special{fp}%
\special{pa 2438 1430}%
\special{pa 2338 1280}%
\special{fp}%
%
\special{pn 8}%
\special{pa 1938 2230}%
\special{pa 2018 2080}%
\special{fp}%
%
\special{pn 8}%
\special{ar 1628 1920 10 10  0.0000000 6.2831853}%
%
\special{pn 8}%
\special{ar 1938 1910 10 10  0.0000000 6.2831853}%
%
\special{pn 8}%
\special{ar 2018 2080 10 10  0.0000000 6.2831853}%
%
\special{pn 8}%
\special{ar 2518 2080 10 10  0.0000000 6.2831853}%
%
\special{pn 8}%
\special{ar 2864 2080 10 10  0.0000000 6.2831853}%
%
\special{pn 8}%
\special{ar 3008 2080 10 10  0.0000000 6.2831853}%
%
\special{pn 8}%
\special{ar 2518 1280 10 10  0.0000000 6.2831853}%
%
\special{pn 8}%
\special{ar 2338 1280 10 10  0.0000000 6.2831853}%
%
\special{pn 8}%
\special{pa 2528 2060}%
\special{pa 2588 1900}%
\special{fp}%
%
\special{pn 8}%
\special{ar 2588 1896 10 10  0.0000000 6.2831853}%
\put(4.8700,-6.4000){\makebox(0,0){$n=1$}}%
\put(4.8700,-14.4000){\makebox(0,0){$n=2$}}%
\put(4.8700,-22.4000){\makebox(0,0){$n=3$}}%
%
\special{pn 8}%
\special{ar 1950 500 0 0  3.9269908 4.7123890}%
%
\special{pn 8}%
\special{sh 1}%
\special{ar 1940 470 10 10 0  6.28318530717959E+0000}%
\end{picture}%
\end{center}
\caption{}
\end{figure}

It is well-known that the ordinary Catalan numbers 
enumerate the number of the rooted trees.
In this section we consider trees associated with graph  $G$.
  
Let $G = (V,E)$ be a directed graph and $G^* =(V^*,E^*)$ the transposed graph of $G$.
A $G$-{\it rooted tree}\, ${\cal T}$ with $n$ edges is an oriented rooted tree with $n$ edges
satisfying the following conditions:
\begin{enumerate}
\item each edge with vertices  is labeled by edges with vertices  of $G^*$,
\item an edge $e^*$ of ${\cal T}$ follows an edge $f^*$ of ${\cal T}$
if and only if   
$e^*$  follows  $f^*$ in the graph $G^*$.
\end{enumerate}
Let $T_n^G$ be the set of all $G$-rooted trees with $n$ edges.
We set 
$$
t_0^G = N, \qquad t_n^G = | T_n^G | \quad \text{ for }  n=1,2,\dots.
$$
For a vertex $u \in V$ and an edge $e \in E$,  
let $T_n^G(u)$ be the set of all $G$-rooted trees whose root is the vertex $u$, and $T_n^G[e]$ the set of all $G$-rooted trees whose root is the source of 
 $e$.  
 Put
\begin{align*}
t_0^G(u) =1,&  \qquad t_n^G(u) = | T_n^G(u) | \quad \text{ for } n=1,2,\dots,\\
t_0^G[e] =1,&  \qquad t_n^G[e] = | T_n^G[e] | \quad \text{ for } n=1,2,\dots.
\end{align*}
\begin{prop}
For $n=0,1,\dots$,
we have  
$t_n^G[e] = d_n^G[e]$ for $e \in E$
and 
$t_n^G(u)= d_n^G(u)$ for $u \in V$,
so that 
$$
t_n^G = d_n^G. 
$$
\end{prop}
\begin{pf}
For a $G$-Dyck path $\gamma$ in $D_n^G$, by considering vectors $(1,1)$ 
in $\gamma$, one gets a $G$-rooted tree.
This correspondence yields bijective mappings 
between $D_n^G$ and $T_n^G$, 
between $D_n^G(u)$ and $T_n^G(u)$, 
and between $D_n^G[e]$ and $T_n^G[e]$.
\end{pf}

\section{Generating functions}
We will next study the generating functions of the sequance
$c_n^G, n=0,1,\dots$.
Let 
$f^G(x) $ 
be the generating function for the sequence
$c_n^G, n=0,1,\dots$, that is defined by
\begin{equation}
f^G(x) = \sum_{n=0}^{\infty} c_n^G x^n  \label{eqn:fg}
\end{equation}
as a formal power series. For a vertex $u,v \in V$, we denote by 
$A_G(v,u)$ the number of edges from $v$ to $u$ in $G$.
By proposition 3.7 the following proposition holds.
\begin{prop} For $n=0,1,\dots$, we have
$$
c_{n+1}^G(u) =\sum_{k=0}^n c_{n-k}^G(u) 
\sum
\begin{Sb}
v \in V \\ 
\end{Sb}
A_G(v,u) c_k^G(v).
$$
\end{prop}
To study the sequence  
$ c_n^G, n=0,1,\dots$
and its generating function
$f^G(x)$,
we provide the generating functions
 for the sequences
$c_n^G(u), n=0,1,\dots$ 
for $u \in V$.
Let $\{ v_1,\dots,v_N \}$ be the vertex set $V$ of $G$.
We put
$$
c_n^G(i)= c_n^G(v_i),\quad
A_G(i,j)  = A_G(v_i,v_j)\quad \text{ for }i,j=1,\dots,N.
$$
Let 
$f_i^G(x) $ 
be the generating function for the sequence
$c_n^G(i), n=0,1,\dots$.
It is defined by
\begin{equation}
f_i^G(x) = \sum_{n=0}^{\infty} c_n^G(i) x^n,\qquad i=1,\dots,N \label{eqn:fig}
\end{equation}
as a formal power series.
The preceding proposition implies the following equalities
\begin{equation}
c_{n+1}^G(i) = \sum_{k=0}^n c_{n-k}^G(i) \sum_{j=1}^N A_G(j,i)c_k^G(j),
\qquad i=1,\dots,N
\end{equation}
so that we have
\begin{prop}
\hspace{6cm}
\begin{enumerate}
\renewcommand{\labelenumi}{(\roman{enumi})}
\item
$
f^G(x) = \sum_{i=1}^N f_i^G(x),
$
\item
\begin{equation}
f_i^G(x)  = 1 + x f_i^G(x) \sum_{j=1}^N A_G(j,i) f_j^G(x), 
\quad \text{ for }\quad i=1,\dots,N.
\end{equation}
\end{enumerate}
\end{prop}
\begin{pf}
(i) The equality is clear.
(ii)
By (6.3),
one has
\begin{align*}
f_i^G(x) - 1  
& = x \sum_{n=0}^{\infty}
c_{n+1}^G(i) x^n \\
& = x \sum_{n=0}^{\infty} 
\sum_{k=0}^n c_{n-k}^G(i) x^{n-k} \sum_{j=1}^N A_G(j,i)c_k^G(j) x^k \\
& = x \sum_{j=1}^N 
\sum_{n=0}^{\infty} A_G(j,i) \sum_{k=0}^n c_{n-k}^G(i) x^{n-k} c_k^G(j) x^k \\
& =  x \sum_{j=1}^N f_i^G(x)  A_G(j,i) f_j^G(x).
\end{align*}
\end{pf}
\begin{lem}
\begin{equation*}
c_n^G (i) \le \| A_G \|_1^n c_n, \qquad i=1,\dots,N,\quad n=0,1,\dots
\end{equation*}
where $\| A_G \|_1 = \max_{1 \le i \le N} \sum_{j=1}^N A_G(j,i)$
and
$c_n = \frac{1}{n+1}\binom{2n}{n}$ the Catalan number.
\end{lem}
\begin{pf}
We will prove the above identity by induction.
Fix $i=1,\dots,N.$
For $n=0$, the inequality is trivial.
For $n=1$, one sees that
$c_1^G(i) = \sum_{j=1}^N A_G(j,i) \le \| A_G \|_1.$
Assume that the inequality holds for all $n\le k$.
As $\sum_{l=0}^k c_{k-l} c_l = c_{k+1}$,
it follows that
\begin{align*}
c_{k+1}^G(i)
& = \sum_{l=0}^k c_{k-l}^G (i)\sum_{j=1}^N A_G(j,i) c_l^G(j) \\
& \le \sum_{l=0}^k \| A_G \|_1^{k-l} c_{k-l} \sum_{j=1}^N A_G(j,i) \| A_G \|_1^l c_l \\ 
& = \| A_G \|_1^{k+1} c_{k+1}  
\end{align*}
Hence the desired inequality holds.
\end{pf}
We denote by $R^G_i$ the radius 
$\frac{1}{\limsup_{n \to \infty}\sqrt[n]{c_n^G(i)}}
$
of convergence of $f_i^G(x)$.
\begin{lem}
Suppose that $G$ is irreducible.
\hspace{3cm}
\begin{enumerate}
\renewcommand{\labelenumi}{(\roman{enumi})}
\item
$R^G_i = R^G_j,\quad i,j=1,\dots,N$.
\item
$\frac{1}{4 \| A_G \|_1} \le R^G_i,\quad i=1,\dots N.$
\end{enumerate}
\end{lem}
\begin{pf}
Put $\alpha_i = \limsup_{n \to \infty}\sqrt[n]{c_n^G(i)}.$

(i)
By the relation (6.3),
one has
$
c_{n+1}^G(i) \ge c_0^G(i) A_G(j,i)c_n^G(j). 
$
Assume that $A_G(j,i) \ne 0$. 
As $c_0^G(i) = 1$, 
one has 
$
c_{n+1}^G(i) \ge c_n^G(j) 
$
so that 
$\alpha_i \ge \alpha_j$.
Since $G$ is irreducible,
one sees that 
$\alpha_i = \alpha_j$
for all $i,j= 1,\dots,N$.

(ii)
It is well-known that
$\limsup_{n \to \infty}\sqrt[n]{c_n} =4$.
By the preceding lemma,
the inequality 
$\alpha_i \le  4 \| A_G \|_1$
is immediate.
\end{pf}
We note that the value $\frac{1}{4 \| A_G \|_1}$
is not best possible in general (see Section 8).

Therefore the functions $f_i^G(x)$ defined by (6.2) 
exists in a neighborhood of the origin, and they satisfy the relations (6.4).
Conversely, 
the following proposition states that the functions are uniquely determined by only the relations (6.4).    
\begin{thm}
A family $f_i(x), i=1,\dots,N$ of functions satisfying  the relations
\begin{equation}
f_i(x)  = 1 + x f_i(x) \sum_{j=1}^N A_G(j,i) f_j(x), 
\qquad  i=1,\dots,N
\label{eqn:f}
\end{equation}
 uniquely exists in a 
neighborhood of the origin,  and they are differentiable.
\end{thm}
\begin{pf}
We first note that by the above equalities one has
$$ 
f_i(0) = 1, \qquad i=1,\dots,N.
$$
Consider the family of polynomials defined by 
\begin{equation*}
F_i(x,y_1,\dots,y_N) = x y_i \sum_{j=1}^N A_G(j,i)y_j - y_i +1, 
\qquad i=1,\dots,N.
\end{equation*}
Put ${\bold y} =(y_1, \dots,y_N) \in {\Bbb R}^N$. 
We set the ${\Bbb R}^N$-valued $C^\infty$-function ${\Bbb F}$ on 
${\Bbb R} \times {\Bbb R}^N$
$$
{\Bbb F}(x,{\bold y})
=
\begin{bmatrix}
F_1(x,y_1,\dots,y_N)\\
F_2(x,y_1,\dots,y_N)\\
\vdots \\
F_N(x,y_1,\dots,y_N)
\end{bmatrix}.
$$
We note that 
${\Bbb F}(0,1,\dots,1) = \bold{0}$.
Since
$$
\frac{\partial F_i}{\partial y_j}=
\begin{cases}
x y_i A(i,i) + x \sum_{k=1}^N A(k,i)y_k -1 & \text{ if } j=i,\\
x A(j,i)y_i & \text{ if } j \ne i,
\end{cases}
$$
one has 
the Jacobian matrix of $\Bbb F$ as
\begin{align*}
\frac{\partial{\Bbb F}}{\partial {\bold y}}
 = & 
\begin{bmatrix}
\frac{\partial F_1}{\partial y_1} & \cdots & \frac{\partial F_1}{\partial y_N}\\\vdots & & \vdots \\
\frac{\partial F_N}{\partial y_1} & \cdots & \frac{\partial F_N}{\partial y_N}  \end{bmatrix}\\
 = & x 
\begin{bmatrix}
y_1 &        & \\
    & \ddots & \\
    &        & y_N
\end{bmatrix}
\begin{bmatrix}
{}^tA(1,1) & \cdots  & {}^tA(1,N)\\
 \vdots    &         & \vdots \\
{}^tA(N,1) & \cdots  & {}^tA(N,N)
\end{bmatrix} \\
+ & 
x 
\begin{bmatrix}
\sum_{k=1}^N {}^tA(1,k)y_k &        & \\
                           & \ddots & \\
                           &        & \sum_{k=1}^N {}^tA(N,k)y_k
\end{bmatrix}
-
\begin{bmatrix}
1 &        & \\
    & \ddots & \\
    &        & 1
\end{bmatrix}.
\end{align*}
Hence we have
$$
\frac{\partial{\Bbb F}}{\partial {\bold y}}(0,1,\dots,1) 
= -
 \begin{bmatrix}
1 &        & \\
    & \ddots & \\
    &        & 1
\end{bmatrix}
\ne 0.
$$
By the implicit function theorem,
one sees the assertion.
\end{pf}
We note that the functions 
$f_1(x),\dots,f_N(x)$
are holomorphic in a neighborhood of the origin. 
 
 Let
$f_1(x),\dots,f_N(x)$
be a family of functions satisfying the equalities (6.5).
They are uniquely defined by 
a neighborhood of the origin by Proposition 6.5. 
Let $f_i^{(n)}(x)$ be the $n$-th derivative of $f_i$.
Since the family of the functions is unique,
the $i$-th $G$-Catalan numbers $c_n^G(i)$
are given by 
$$
 c_n^G(i)= \frac{f_i^{(n)}(0)}{n!},
 \qquad n=0,1,\dots.
 $$
We henceforth assume that $G$ is irreducible.
 We denote by $R_G$ the radius $R_i^G$ of convergence of $f_i^G(x)$ 
 as in Lemma 6.4. 
 We put
 $I_{R_G} = \{ x \in {\Bbb R} \mid | x | < R_G \}.$
 \begin{lem}
 \hspace{3cm}
\begin{enumerate}
\renewcommand{\labelenumi}{(\roman{enumi})}
\item $f_i(x) \ne 0$ for $ x \in I_{R_G}$ 
and $ f_i(x) > 1$ for $0 < x \in I_{R_G}$.
\item There exists $M > 0$ such that 
$| f_i(x) | < M$ for all $ x \in I_{R_G}, \, i=1,\dots,N$.
\end{enumerate}
 \end{lem}
 \begin{pf}
(i) Since $f_i(x) = f^G_i(x)$ on $x \in I_{R_G}$, 
the assertion (i) is clear by definition of $f_i^G(x)$. 

(ii)
Suppose that there exists $j=1,\dots,N$ and $x_n \in I_{R_G}$
such that $\lim_{n \to \infty} |f_j (x_n) | = \infty.$
Since
$| f_j(x_n) | \le f_j (| x_n | )$, 
we may assume that $x_n > 0$ by considering $ |x_n | $ instead of $x_n$.
Take $i = 1,\dots,N$ such that 
$A_G(j,i) = 1$. 
As we have
$$
f_i(x_n) -1 \ge x_n f_i(x_n) f_j(x_n),
$$ 
the inequality 
$$
1 \ge 1 - \frac{1}{f_i(x_n)} \ge x_n f_j(x_n) \ge 0
$$
holds by (i).
By hypothesis, 
one sees that 
$\lim_{n \to \infty} x_n = 0$,
a contradiction to the fact
$ 1 = f_j(0) $ with the continuity of $f_j$ at $0$.
 \end{pf}
\begin{prop}
The functions $f_i(x),i=1,\dots,N$ satisfying \eqref{eqn:f}
can be defined at $x = R_G$, and they are lower semi-continuous at $x= R_G$.
 \end{prop}
\begin{pf}
Take an increasing sequence $x_n$ in $I_{R_G}$ such that 
$x_n \uparrow R_G$.
For each $i=1,\dots,N$, 
the sequence $\{ f_i (x_n) \}_{n=0,1, \dots }$ is increasing and bounded by 
the preceding lemma so that 
 the functions $f_i(x)$ 
can be defined at $x = R_G$, and they are lower-continuous at $x= R_G$.
\end{pf}
Therefore we have 
\begin{thm}
The family $f_i(x),i=1,\dots, N $ of functions
satisfying \eqref{eqn:f}
uniquely exists on $[-R_G, R_G]$ for some $R_G >\frac{1}{4 \|A \|_1}$.  
They can not be defined outside of the interval $[-R_G, R_G]$.
\end{thm}
We note the following proposition.
\begin{prop}
If 
$
\{ j \mid A(j,i_1) =1 \} =\{ j \mid A(j,i_2) =1 \},
$
then we have
$f_{i_1}(x) = f_{i_2}(x).$
\end{prop}
\begin{pf}
Put
$
g_i(x) = \sum_{j=1}^N A_G(j,i)f_j(x). 
$
As $f_i(x)$ is equal to
$\frac{1}{1 - xg_i(x)}$ for $x \ne 0$,
one sees 
$f_{i_1}(x) = f_{i_2}(x)$
for $x \ne 0$ by hypothesis.
As
$f_{i_1}(0) = f_{i_2}(0)=1$,
we have
$f_{i_1}(x) = f_{i_2}(x).$
\end{pf}

\section{The radius of convergence of the generating functions}

In this section, we will study how to find the radius $R_G$ of 
convergence of the functions $f_i(x), i=1,\dots,N$ 
satisfying \eqref{eqn:f}.
For an $N \times N$  matrix $A$ and $t = (t_i)_{i=1}^N \in {\Bbb R}^N$,
we set 
$$
x(t)_i = t_i - \sum_{k=1}^N t_k A(k,i) t_i,
\qquad
(tAt)(i,j) = t_i A(i,j)t_j.
$$
for $i,j=1,\dots,N$.
Hence we have an $N\times N$ matrix $tAt$. We set 
$$
C_A = \{ t=(t_i)_{i=1}^N \in {\Bbb R}^N \mid
 t_i > 0, \, x(t)_i = x(t)_j > 0 \text{ for }
 i,j=1,\dots,N \}.
 $$
We note that 
$$
C_{A_G} = \{ (x f_i^G(x))_{i=1}^N \in {\Bbb R}^N \mid
 0 < x \in I_{R_G}  \}
 $$
 by the relations (6.4) and Theorem 6.8.  
Let $\widetilde{A}$
be the $N \times N$ matrix defined by
$$
\widetilde{A}(i,j) = L_A -A(i,j),\qquad i,j = 1,\dots,N
$$
where $L_A = \max_{i,j}A(i,j)$.
For 
$t=(t_i)_{i=1}^N \in C_A$,
the value $x(t)_i$ does not depend on $i=1,\dots,n$.
We denote it by
$x(t)$.
We say that $A$ satisfies {\it condition}\, (C)
if 
$$
\ker(t A t - x(t) ) \cap \ker ( t \widetilde{A} t + x(t) ) = \{ 0 \}
$$
for all $t=(t_i)_{i=1}^N \in C_A$.
\begin{lem}
A matrix $A$ satisfies condition (C) 
if and only if
$$
\ker(t A t - x(t) ) \cap {\bold 1}^{\perp} = \{ 0 \}.
$$
for all
$t=(t_i)_{i=1}^N \in C_A$,
where
$
{\bold 1}^{\perp} 
= \{ (r_i )_{i=1}^N \in {\Bbb R}^N \mid \sum_{i=1}^N r_i= 0 \}.
$
\end{lem}
\begin{pf}
We will show that for $(r_i)_{i=1}^N \in \ker(t A t - x(t) )$,
the vector $(r_i)_{i=1}^N$ belongs to $ \ker ( t \widetilde{A} t + x(t) ) $
if and only if 
$\sum_{i=1}^N r_i= 0 $.
As $(r_i)_{i=1}^N \in \ker(t A t - x(t) )$,
one has
$$
\sum_{i=1}^N x(t) r_i 
 = \sum_{i=1}^N t_i r_i - \sum_{i,k=1}^N t_k A(k,i)t_i r_i
  = \sum_{i=1}^N t_i r_i - \sum_{i=1}^N x(t) r_i 
$$
so that
\begin{equation*}
2 x(t) \sum_{i=1}^N r_i  = \sum_{i=1}^N t_i r_i. 
\end{equation*}
Since $x(t) >0$, one knows that 
$\sum_{i=1}^N r_i= 0 $
if and only if 
$\sum_{i=1}^N t_i r_i= 0$.
Now 
$$
t \widetilde{A} t 
\begin{bmatrix}
r_1\\
\vdots \\
r_N \\
\end{bmatrix}
= 
\begin{bmatrix}
L_A t_1\sum_{j=1}^N t_j r_j - t_1 \sum_{j=1}^N A(1,j)t_j r_j \\
\vdots \\
L_A t_N \sum_{j=1}^N t_j r_j - t_N \sum_{j=1}^N A(N,j)t_j r_j \\
\end{bmatrix}
= 
L_A \sum_{j=1}^N t_j r_j
\begin{bmatrix}
 t_1  \\
\vdots \\
 t_N  \\
\end{bmatrix}
- x(t)
\begin{bmatrix}
 r_1 \\
\vdots \\
r_N \\
\end{bmatrix}
$$
Hence
$\sum_{i=1}^N r_i= 0 $
if and only if 
$(r_i)_{i=1}^N$ belongs to $ \ker ( t \widetilde{A} t + x(t) )$.
\end{pf}
It is easy to see that
the matrices
$$
[N], \quad
\begin{bmatrix} 
1      & \cdots & 1 \\
\vdots &        & \vdots \\ 
1      & \cdots & 1
\end{bmatrix},
\quad
\begin{bmatrix} 
1      &  1 \\
1      & 0 
\end{bmatrix}
$$
satisfy condition (C).

We will prove the following theorem:
\begin{thm}\label{th:1}
Suppose that a matrix $A_G$ satisfies condition $(C)$. 
Let $f_i(x), i=1,\dots,N$ be the functions satisfying \eqref{eqn:f}.
If a real number $x_0 \in \Bbb R$ 
is the radius $R_G$ of convergence of the functions 
$f_i(x)$, then 
there exist positive real numbers $t_1,\dots,t_N$
such that
\begin{equation*}
\det (t A_G t - x_0 ) = 0, \qquad
x_0 = t_i - \sum_{j=1}^N t_j A_G(j,i) t_i
\quad \text{ for } i = 1,\dots, N
\end{equation*} 
where
$t A_G t $ is the $N \times N$ matrix defined by 
$ t A_G t = [ t_i A_G(i,j) t_j ]_{i,j = 1,\dots N}.$
In particular, there exists
an eigenvector $[s_i]_{i=1,\dots,N}$
of the matrix $t A_G t$ for the eigenvalue $x_0$ with
$\sum_{i=1}^N s_i = 1$ such that
$$
x_0 = \frac{1}{2} \sum_{i=1}^N t_i s_i.
$$
\end{thm}
Therefore the radius of convergence of $f_i(x)$ is algebraically determined 
as a solution of an eigenvalue problem for the matrix $ tA_G t$ with some conditions.
\begin{pf}
Put $t_i = x f_i(x), i=1,\dots,N$
in \eqref{eqn:f} so that we have equalities 
\begin{equation}
x = t_i - t_i \sum_{j=1}^N A_G(j,i) t_j, \qquad i=1,\dots, N.
\label{eqn:ti}
\end{equation}
This implies that 
$(t_i)_{i=1}^N$ belongs to $C_A$ for $ 0 < x \in I_{R_G}$.
Consider $x=x(t_1,\dots,t_N)$ as a function of $(t_1,\dots,t_N)$
so that the radius $R_G$ 
of convergence is the maximum value of the function 
$x=x(t_1,\dots,t_N)$.
 Put 
$$
 \psi_i(t_1,\dots,t_N) 
 = t_i - t_i \sum_{j=1}^N A_G(j,i) t_j, \qquad i=1,\dots, N
$$
and
\begin{align*}
f(t_1,\dots,t_N) & = \psi_1(t_1,\dots,t_N), \\
g_2(t_1,\dots,t_N) & = \psi_1(t_1,\dots,t_N)-\psi_2(t_1,\dots,t_N), \\
                   & \cdots                                        \\
g_N(t_1,\dots,t_N) & = \psi_1(t_1,\dots,t_N)-\psi_N(t_1,\dots,t_N). 
\end{align*}
Hence the radius $R_G$ is obtained by solving the constrained extremal problem
of $f$ with constrained conditions:
$$
g_2(t_1,\dots,t_N)= \cdots = g_N(t_1,\dots,t_N) =0.
$$
Suppose that 
$$
\sum_{j=2}^N c_j \frac{\partial g_j}{\partial t_i}(t_1,\dots,t_N) = 0,
\qquad
i=1,\dots,N
$$
for some $c_j \in {\Bbb R}, j=2,\dots,N$.
One then has
\begin{equation*}
\begin{bmatrix}
\frac{\partial \psi_1}{\partial t_1}(t_1,\dots,t_N) & 
\frac{\partial \psi_2}{\partial t_1}(t_1,\dots,t_N) &
\cdots &
\frac{\partial \psi_N}{\partial t_1}(t_1,\dots,t_N) \\
\frac{\partial \psi_1}{\partial t_2}(t_1,\dots,t_N) & 
\frac{\partial \psi_2}{\partial t_2}(t_1,\dots,t_N) &
\cdots &
\frac{\partial \psi_N}{\partial t_2}(t_1,\dots,t_N) \\
\vdots & 
\vdots &
\cdots &
\vdots \\
\frac{\partial \psi_1}{\partial t_N}(t_1,\dots,t_N) & 
\frac{\partial \psi_2}{\partial t_N}(t_1,\dots,t_N) &
\cdots &
\frac{\partial \psi_N}{\partial t_N}(t_1,\dots,t_N) \\
\end{bmatrix}
\begin{bmatrix}
\sum_{j=2}^N c_j \\
 - c_2 \\
 \vdots \\
 - c_N
\end{bmatrix}
=0.
\end{equation*}
As one sees 
\begin{equation*}
\frac{\partial \psi_i}{\partial t_j} =
\begin{cases}
1 - \sum_{k=1}^N A_G(k,i) t_k - A_G(i,i) t_i & \text{ if } i=j \\
- A_G(j,i) t_i & \text{ if } i \ne j,
\end{cases}
\end{equation*}
by putting 
$$
r_1 = \sum_{j=2}^N c_j,\quad 
r_2=  - c_2,\quad \dots,\quad
r_N = - c_N,
$$
one has
$$
(1 - \sum_{k=1}^N t_k A_G(k,i)) r_i= \sum_{j=1}^N  A_G(i,j) t_j r_j, 
\qquad i=1,\dots,N.
$$
By \eqref{eqn:ti}, one has
$$
  x
\begin{bmatrix}
r_1 \\
\vdots \\
r_N
\end{bmatrix}
=  
t A_G t
\begin{bmatrix}
r_1 \\
\vdots \\
r_N
\end{bmatrix},
\qquad
\sum_{i=1}^N r_i = 0.
$$
Since the matrix $A_G$ satisfies condition $(C)$,
one gets $r_i=0, i=1,\dots,N$.
This means that 
the rank of the matrix 
$
[\frac{\partial g_j}{\partial t_i}(t_1,\dots,t_N)]_{i=1,\dots,N, \,j=2,\dots,N}
$ is $N-1$.
Now suppose that the function $f(x)$ takes its maximum value $x_0$ at 
$(t_1,\dots,t_N)$ under the conditions that 
$
g_i(t_1,\dots,t_N)  = 0
$
for
$
i=2,\dots,N.
$
Then there exists real numbers
$\lambda_1,\dots,\lambda_N$  
such that
$$
\frac{\partial f}{\partial t_i} (t_1,\dots,t_N)
 + \sum_{j=2}^N \lambda_j \frac{\partial g_j}{\partial t_i}(t_1,\dots,t_N)
  = 0, 
\qquad i=1,\dots,N.\label{eqn:pgi}
$$
It then follows that
\begin{equation*}
\frac{\partial \psi_1}{\partial t_i}(t_1,\dots,t_N)
 + \sum_{j=2}^N 
  \lambda_j ( \frac{\partial \psi_1}{\partial t_i}(t_1,\dots,t_N)
  -  \frac{\partial \psi_j}{\partial t_i}(t_1,\dots,t_N)) = 0
\end{equation*}
so that 
\begin{equation*}
\begin{bmatrix}
\frac{\partial \psi_1}{\partial t_1}(t_1,\dots,t_N) & 
\frac{\partial \psi_2}{\partial t_1}(t_1,\dots,t_N) &
\cdots &
\frac{\partial \psi_N}{\partial t_1}(t_1,\dots,t_N) \\
\frac{\partial \psi_1}{\partial t_2}(t_1,\dots,t_N) & 
\frac{\partial \psi_2}{\partial t_2}(t_1,\dots,t_N) &
\cdots &
\frac{\partial \psi_N}{\partial t_2}(t_1,\dots,t_N) \\
\vdots & 
\vdots &
\cdots &
\vdots \\
\frac{\partial \psi_1}{\partial t_N}(t_1,\dots,t_N) & 
\frac{\partial \psi_2}{\partial t_N}(t_1,\dots,t_N) &
\cdots &
\frac{\partial \psi_N}{\partial t_N}(t_1,\dots,t_N) \\
\end{bmatrix}
\begin{bmatrix}
1 + \sum_{j=2}^N \lambda_j \\
 - \lambda_2 \\
 \vdots \\
 - \lambda_N
\end{bmatrix}
=0.
\end{equation*}
Put
$$
s_1 = 1 + \sum_{j=2}^N \lambda_j, \quad  
s_2 = - \lambda_2, \quad \dots, \quad s_N = - \lambda_N.
$$
Similarly to the above discussions,
we have
$$
  x_0
\begin{bmatrix}
s_1 \\
\vdots \\
s_N
\end{bmatrix}
=  
t A_G t
\begin{bmatrix}
s_1 \\
\vdots \\
s_N
\end{bmatrix}, \qquad
\sum_{j=1}^N s_j =1.
$$
In this case we have
$$
x_0 = \sum_{i=1}^N t_i s_i - \sum_{i,j =1}^N t_j A_G(j,i) t_i s_i 
    = \sum_{i=1}^N t_i s_i - x_0.
$$
so that we get 
$$
x_0 = \frac{1}{2} \sum_{i=1}^N t_i s_i.
$$
Therefore we obtain the assertion.
\end{pf}
\section{Integral formulae of the $G$-Catalan numbers}
The $G$-Catalan numbers $c_n^G(i)$ are coefficients of its generating functions $f_i^G(x)$ so that one has
$$
c_n^G(i) = \frac{1}{n!} \frac{d^{n} f_i^G (x)}{dx^n}(0) 
=\frac{1}{2\pi \sqrt{-1}} \int_C \! \! \frac{f_i^G(z)}{z^{n+1}} dz
$$
where $C$ denotes a positively oriented closed curve around the origin
in the complex plane.  
In this section, we present a formula of  $c_n^G(i)$
by using the relations (6.4).
We set 
\begin{align*}
F_i(w_1,\dots,w_N) & = (w_i +1) \sum_{j=1}^N A_G(j,i) (w_j +1),\\
F_i^n(w_1,\dots,w_N) & = F_i(w_1,\dots,w_N)^n,\qquad i=1,\dots,N.
\end{align*}
We will prove the following integral formulae 
\begin{thm}
\begin{equation}
c_n^G(i) =  \frac{1}{2\pi \sqrt{-1} n} \int_C \! \! 
\frac{F_j^n(w_1,\dots,w_N)}{w_j^n} dw_i,\qquad i,j =1,\dots,N,
\end{equation}
where the above integral is a contour integral along a positively oriented closed curve $C$
around the origin in the complex plane.
\end{thm}
\begin{pf}
Put
$ w_i(x) = f_i^G(x) -1$.
The relations (6.4) go to
$$
w_i(x) = x \sum_{j=1}^N (w_j(x) + 1 )A_G(j,i) (w_i(x) + 1) 
= x F_i(w_1(x),\dots, w_N(x)).
$$
As
$$
\lim_{x \to 0} \frac{w_i(x)}{x} = F_i(w_1(0),\dots, w_N(0)) =
\sum_{j=1}^N A_G(j,i) > 0,
$$
the rotation number of $w_i(x)$ is $1$.
One sees that
$$
w'_i (z)  = f'_i(z)
$$
and
$$
\frac{1}{z} =  \frac{ F_j(w_1(z),\dots, w_N(z))}{w_j(z)}.
$$
It then folows that
\begin{align*}
c_n^G(i) 
& =  \frac{1}{2\pi \sqrt{-1} n} \int_C \! \!
\frac{f'_i(z)}{z^{n}} dz \\
& =  \frac{1}{2\pi \sqrt{-1} n} \int_C \! \!
\frac{ F_j^n(w_1(z),\dots, w_N(z))}{w_j^n(z)} w'_i (z) dz.
\end{align*} 
\end{pf}

\section{Catalan numbers associated to KMS states}
In this section, we enumerate $G$-Catalan numbers by using KMS states on the Cuntz-Krieger algebra ${\cal O}_{A^G}$ for the gauge actions.
We denote by ${\Bbb T}$ the group of complex numbers with modulus one.
The gauge action $\alpha^G$ is an action of ${\Bbb T}$ to  
the automorphisms on ${\cal O}_{A^G}$ defined by 
$\alpha^G_z(S_e)= z S_e$ for $z \in {\Bbb T},e \in E$.
For a real number $\beta\in {\Bbb R}$, 
a state $\varphi$ on ${\cal O}_{A^G}$ is called a KMS state at inverse temperature $\beta$ if the following equality holds
$$
\varphi(ab) = \varphi(b \alpha^G_{t + i\beta}(a)),\qquad t \in {\Bbb R} 
$$
for $a $ in the  dense analytic elements of 
${\cal O}_{A^G}$ and $b \in {\cal O}_{A^G}$.  
In \cite{EFW84}, it has been proved that 
under the condition that the matrix $A^G$ is aperiodic,
KMS state exists if and only if $\beta$ is $\log r_G$, 
and the admitted KMS state is unique, 
where $r_G$ is the Perron-Frobenius eigenvalue for the matrix $A^G$.
In what follows, we fix a KMS state $\varphi$ on ${\cal O}_{A^G}$.
We define the $(G,\varphi)$-Catalan numbers
$c_n^{G,\varphi}, n=0,1,\dots $ by setting
$$
c_0^{G,\varphi} = 1,\qquad  
c_n^{G,\varphi} = \sum_{X \in B_n^G} \varphi(\pi_G(X)), \quad n=1,2,\dots.
$$
These numbers are not necessarily integers.
For $i=1,\dots,N$, we put the
$(G,\varphi)$-Catalan numbers $c_n^{G,\varphi}(i), n=0,1,\dots $ 
rooted at the vertex $v_i$ by setting
$$
c_0^{G,\varphi}(i) = \varphi(P_{v_i}),\qquad  
c_n^{G,\varphi}(i) = \sum_{X \in B_n^G(v_i)} \varphi(\pi_G(X)), 
\quad n=1,2,\dots.
$$
By the preceding discussions, one  has
\begin{lem}
For $n=0,1,\dots,$ and $i=1,\dots,N$, we have
\hspace{6cm}
\begin{enumerate}
\renewcommand{\labelenumi}{(\roman{enumi})}
\item
$c_n^{G,\varphi} = \sum\limits_{i=1}^N c_n^{G,\varphi}(i).$
\item
$c_n^{G,\varphi}(i) = c_n^G(i) \varphi(P_{v_i}).$
\item
$
c_{n+1}^{G,\varphi}(i) =\sum\limits_{k=0}^n  
\sum\limits_{j=1}^N c_k^G(j)
A_G(j,i)c_{n-k}^{G,\varphi}(i).
$
\end{enumerate}
\end{lem}
\begin{pf}
(i) is clear.

(ii) For $X \in B^G_n$, the word $X$ belongs to $B_n^G(i)$ if and only if
$v(X) = v_i$. 
The latter condition 
is equivalent to the condition $\pi_G(X) = P_{v_i}$.
This implies the assertion (ii).

(iii) The desired equality comes from (ii) with (6.3).
\end{pf}
  The generating functions $f^{G,\varphi}(x)$ and 
$f^{G,\varphi}_i(x)$ for the sequences 
$c_n^{G,\varphi}$ and $c_n^{G,\varphi}(i)$ respectively
 are also defined by 
$$
f^{G,\varphi}(x)  = \sum_{n=0}^\infty c_n^{G,\varphi} x^n, \qquad
f_i^{G,\varphi}(x) = \sum_{n=0}^\infty c_n^{G,\varphi}(i) x^n 
\qquad \text{ for } i=1,\dots,N.
$$
Then the following lemma is direct from the preceding lemma.
\begin{lem}
For $i=1,\dots,N$, we have
\hspace{6cm}
\begin{enumerate}
\renewcommand{\labelenumi}{(\roman{enumi})}
\item
$f^{G,\varphi}(x) = \sum\limits_{j=1}^N f^{G,\varphi}_j(x).$
\item
$f^{G,\varphi}_i(x) = f^G_i(x) \varphi(P_{v_i}).$
\item
$
f^{G,\varphi}_i(x) = \varphi(P_{v_i}) +  
x  \sum\limits_{j=1}^N f^G_j(x) A_G(j,i)f^{G,\varphi}_i(x).
$
\end{enumerate}
\end{lem}
The radii of convergence of the functions 
$f^{G,\varphi}(x), f^{G,\varphi}_i(x)$
coincide with that of $f^G(x)$ if $G$ is irreducible.

Therefore we have
\begin{prop}
Suppose that $A^G$ is aperiodic.
Let ${[t_e]}_{e \in E}$ be the positive eigenvector for the Perron-Frobenius eigenvalue $r_G$ of the matrix $[A^G(e,f)]_{e,f\in E}$
satisfying $\sum_{e\in E}t_e =1$.
For a vertex $v_i\in V,i=1,\dots,N$, take an edge $e_i \in E$ such that 
$t(e_i) = v_i$.
Then we have
\begin{align*}
c^{G,\varphi}_n(i) & = r_G c^G_n(i)  t_{e_i},
\qquad  
c^{G,\varphi}_n = r_G \sum_{i=1}^N c^G_n(i) t_{e_i}\\
\intertext{and hence}
f^{G,\varphi}_i(x) & = r_G f^G_i(x) t_{e_i},
\qquad
f^{G,\varphi}(x) = r_G \sum_{i=1}^N f^G_i(x) t_{e_i}.
\end{align*} 
\end{prop}
\begin{pf}
As in \cite{EFW84}, the vector 
$[\varphi(S_e S_e^*)]_{e \in E}$
is the unique positive eigenvector for the Perron-Frobenius eigenvalue 
$r_G$ of the matrix $[A^G(e,f)]_{e,f \in E}$  satisfying 
$\sum_{f \in E} \varphi(S_f S_f^*) =1$. 
Hence we have $\varphi(S_e S_e^*) = t_e$ for $e\in E$.
It follows that
$$
\varphi(S_e^* S_e ) = \sum_{f \in E}A^G(e,f)\varphi(S_f S_f^*)=r_G t_e
$$
and 
$
\varphi(P_{v_i}) =\varphi(S_{e_i}^* S_{e_i}) = r_G t_{e_i}.
$
By the preceding lemma, one has
$
c^{G,\varphi}_n(i) =  r_G c^G_n(i) t_{e_i}
$ 
and 
$ 
c^{G,\varphi}_n = r_G \sum_{i=1}^N c^G_n(i) t_{e_i}.
$
\end{pf}

\section{Examples}

1.
Consider the following graph $G_1$ 
having $N$-loops with a single vertex.

\begin{figure}[htbp]
\begin{center}
\unitlength 0.1in
\begin{picture}(  8.4600, 12.8900)( 18.0900,-14.0100)
%
\special{pn 8}%
\special{ar 2232 1382 14 20  0.0000000 6.2831853}%
%
\special{pn 8}%
\special{ar 2232 1308 40 60  1.7046073 6.2831853}%
\special{ar 2232 1308 40 60  0.0000000 1.2565644}%
%
\special{pn 8}%
\special{pa 2248 1364}%
\special{pa 2244 1364}%
\special{fp}%
\special{sh 1}%
\special{pa 2244 1364}%
\special{pa 2314 1368}%
\special{pa 2296 1352}%
\special{pa 2304 1328}%
\special{pa 2244 1364}%
\special{fp}%
%
\special{pn 8}%
\special{ar 2232 1196 118 176  1.6680495 6.2831853}%
\special{ar 2232 1196 118 176  0.0000000 1.4393276}%
%
\special{pn 8}%
\special{pa 2252 1370}%
\special{pa 2248 1370}%
\special{fp}%
\special{sh 1}%
\special{pa 2248 1370}%
\special{pa 2314 1390}%
\special{pa 2300 1370}%
\special{pa 2314 1350}%
\special{pa 2248 1370}%
\special{fp}%
%
\special{pn 8}%
\special{ar 2232 746 424 634  1.6038812 6.2831853}%
\special{ar 2232 746 424 634  0.0000000 1.5282688}%
%
\special{pn 8}%
\special{pa 2254 1378}%
\special{pa 2250 1378}%
\special{fp}%
\special{sh 1}%
\special{pa 2250 1378}%
\special{pa 2318 1398}%
\special{pa 2304 1378}%
\special{pa 2318 1358}%
\special{pa 2250 1378}%
\special{fp}%
%
\special{pn 8}%
\special{sh 1}%
\special{ar 2232 182 10 10 0  6.28318530717959E+0000}%
\special{sh 1}%
\special{ar 2232 296 10 10 0  6.28318530717959E+0000}%
\special{sh 1}%
\special{ar 2232 408 10 10 0  6.28318530717959E+0000}%
\special{sh 1}%
\special{ar 2232 520 10 10 0  6.28318530717959E+0000}%
\special{sh 1}%
\special{ar 2232 632 10 10 0  6.28318530717959E+0000}%
\special{sh 1}%
\special{ar 2232 746 10 10 0  6.28318530717959E+0000}%
\special{sh 1}%
\special{ar 2236 858 10 10 0  6.28318530717959E+0000}%
\special{sh 1}%
\special{ar 2232 970 10 10 0  6.28318530717959E+0000}%
\special{sh 1}%
\special{ar 2232 970 10 10 0  6.28318530717959E+0000}%
\end{picture}%
\end{center}
\caption{}
\end{figure}

 Since $A_{G_1} = [N]$, we have  
$$
f^{G_1}(x) = f^{G_1}_1(x) 
\quad \text{ and } \quad f^{G_1}_1 (x) -1 = x N f^{G_1}_1(x)^2
$$
so that 
$$
f^{G_1}(x) = \frac{1 -\sqrt{1 - 4N x}}{2Nx}.
$$
We choose the minus sign before the root of $1 - 4N x$ because 
$\lim_{x \to 0} f^{G_1}(x) = 1.$ 
Although  one directly konows that 
the radius $R_{G_1}$ of convergence of $f^{G_1}(x)$
is equal to $\frac{1}{4N}$ because of $1 -4 x N \ge 0$,
we will see $R_{G_1}= \frac{1}{4N}$ by using Theorem 7.2
as in the following way.
Suppose that 
 a positive real number $x_0$ is $R_{G_1}$. 
There exists $t >0$ such that 
 $$ 
 tNt - x_0 = 0, \qquad x_0 = t - tNt.
 $$
 By these equations we have 
 $
 t= \frac{1}{2N}
 $
 and
 $
 x_0 = \frac{1}{4N}.
 $
 We will next compute $c_n^{G_1}$.
 Althogh by the Newton Binomial formula for $(1 -4N x)^{\frac{1}{2}}$
 one may easily compute $c_n^{G_1}$, 
 we will use Theorem 8.1 as follows.
 Put
 $F(w) =  N (w+1)^2$
 so that we have
 $$
 c_n^{G_1} = \frac{1}{2 \pi n \sqrt{-1}} \int_C \! \!
 \frac{F^n(w)}{w^n} dw = \frac{1}{n}N^n \binom{2n}{n-1} = N^n c_n.
 $$    
 Since $\varphi(P_{v_1}) =1$, the equality 
 $c^{G_1,\varphi}_n = c^{G_1}_n$ holds.

2. Let $G_2$ be the directed graph with $N$ vertices such that for any ordered pair of two vertices $u, v$, there uniquely exists an edge from $u$ to $v$.
Hence the transition matrix $A_{G_2}$ is  
$
\begin{bmatrix} 
1      & \cdots & 1 \\
\vdots &        & \vdots \\ 
1      & \cdots & 1
\end{bmatrix}.
$
The generating functions $f^{G_2}_i$ satisfy the following equalities
\begin{equation*}
f^{G_2}_i(x)  -1 = x f^{G_2}_i(x) \sum_{j=1}^N f^{G_2}_j(x), 
\qquad i=1,\dots,N.
\end{equation*}
By Proposition 6.9,
one has
$f^{G_2}_i(x) = f^{G_2}_j(x)$ for $i,j=1,\dots,N$.
By the above equalities,
we have
$$
f^{G_2}_i(x) -1 = x N {f^{G_2}_i}^2(x).
$$
Therefore we have
$$
f^{G_2}_i(x) = f^{G_1}(x) = \frac{1 -\sqrt{1 - 4N x}}{2Nx}.
$$
Hence we have
\begin{align*}
R_{G_2} & = R_{G_1} = \frac{1}{4N},\\
f^{G_2}(x) &= N f^{G_2}_i(x) = \frac{1 -\sqrt{1 - 4N x}}{2x},\\
c_n^{G_2} & = \sum_{i=1}^N c_n^{G_2}(i) = N c_n^{G_1}
=\frac{1}{n}N^{n+1} \binom{2n}{n-1} = N^{n+1} c_n
\end{align*}
By the equality
$\varphi(P_{v_i}) = \frac{1}{N}$, 
we note
$$
c^{(G_2, \varphi)}_n =c^{G_2}_n \cdot \frac{1}{N} = c^{G_1}_n = N^n c_n.
$$

3. Let $G_3$ be the directed graph with 
$
A_{G_3} =
\begin{bmatrix} 
1      &  1 \\
1      & 0
\end{bmatrix}.
$

\begin{figure}[htbp]
\begin{center}
\unitlength 0.1in
\begin{picture}( 21.9100,  6.9000)( 16.5000,-24.0600)
%
\special{pn 8}%
\special{ar 2092 2064 78 66  0.0000000 6.2831853}%
%
\special{pn 8}%
\special{ar 3762 2058 80 66  0.0000000 6.2831853}%
%
\special{pn 8}%
\special{ar 1832 2058 182 154  0.3190696 6.0756891}%
%
\special{pn 8}%
\special{pa 2002 2112}%
\special{pa 2006 2106}%
\special{fp}%
\special{sh 1}%
\special{pa 2006 2106}%
\special{pa 1958 2158}%
\special{pa 1982 2154}%
\special{pa 1994 2176}%
\special{pa 2006 2106}%
\special{fp}%
%
\special{pn 8}%
\special{ar 2946 2128 818 278  6.2480092 6.2831853}%
\special{ar 2946 2128 818 278  0.0000000 3.1562795}%
%
\special{pn 8}%
\special{pa 3762 2122}%
\special{pa 3762 2118}%
\special{fp}%
\special{sh 1}%
\special{pa 3762 2118}%
\special{pa 3742 2186}%
\special{pa 3762 2172}%
\special{pa 3782 2186}%
\special{pa 3762 2118}%
\special{fp}%
%
\special{pn 8}%
\special{pa 2134 2006}%
\special{pa 2136 1974}%
\special{pa 2148 1946}%
\special{pa 2166 1918}%
\special{pa 2188 1896}%
\special{pa 2214 1876}%
\special{pa 2240 1858}%
\special{pa 2268 1842}%
\special{pa 2296 1830}%
\special{pa 2326 1816}%
\special{pa 2356 1804}%
\special{pa 2386 1794}%
\special{pa 2416 1784}%
\special{pa 2448 1776}%
\special{pa 2478 1768}%
\special{pa 2510 1760}%
\special{pa 2542 1754}%
\special{pa 2572 1748}%
\special{pa 2604 1742}%
\special{pa 2636 1738}%
\special{pa 2668 1734}%
\special{pa 2700 1730}%
\special{pa 2732 1728}%
\special{pa 2762 1724}%
\special{pa 2794 1722}%
\special{pa 2826 1720}%
\special{pa 2858 1720}%
\special{pa 2890 1718}%
\special{pa 2922 1716}%
\special{pa 2954 1716}%
\special{pa 2986 1718}%
\special{pa 3018 1718}%
\special{pa 3050 1718}%
\special{pa 3082 1720}%
\special{pa 3114 1722}%
\special{pa 3146 1726}%
\special{pa 3178 1728}%
\special{pa 3210 1730}%
\special{pa 3242 1736}%
\special{pa 3274 1740}%
\special{pa 3306 1744}%
\special{pa 3336 1750}%
\special{pa 3368 1756}%
\special{pa 3400 1764}%
\special{pa 3430 1770}%
\special{pa 3462 1778}%
\special{pa 3492 1788}%
\special{pa 3522 1796}%
\special{pa 3554 1806}%
\special{pa 3582 1820}%
\special{pa 3612 1832}%
\special{pa 3640 1846}%
\special{pa 3668 1862}%
\special{pa 3694 1880}%
\special{pa 3718 1902}%
\special{pa 3740 1926}%
\special{pa 3758 1952}%
\special{pa 3768 1982}%
\special{pa 3768 2000}%
\special{sp}%
%
\special{pn 8}%
\special{pa 2134 1992}%
\special{pa 2134 2006}%
\special{fp}%
\special{sh 1}%
\special{pa 2134 2006}%
\special{pa 2154 1940}%
\special{pa 2134 1954}%
\special{pa 2114 1940}%
\special{pa 2134 2006}%
\special{fp}%
\end{picture}%
\end{center}
\caption{}
\end{figure}
We then have 
$$
f^{G_3}(x) = f^{G_3}_1(x) + f^{G_3}_2(x).
$$
Put
$f_i(x) = f^{G_3}_i(x), i=1,2$.
They satisfy the following relations:
\begin{align*}
f_1(x) -1 & = x (f_1(x) + f_2(x)) f_1(x),\\
f_2(x) -1 & = x f_1(x)  f_2(x)   
\end{align*}
so that the equalities
\begin{equation}
f_2(x)^2  = f_1(x),\qquad
x f_2(x)^3  -f_2(x) + 1  = 0    
\end{equation}
hold.
We will compute the radius $R_{G_3}$ of convergence of the function
$f^{G_3}(x)$.
Suppose that 
 a positive real number $x_0$ is $R_{G_3}$.
 By Theorem 7.2, 
 there exist $t_1, t_2 >0$ such that 
 \begin{gather*}
 \det( 
\begin{bmatrix}
 t_1^2 & t_1 t_2 \\
 t_2 t_1 & 0 
\end{bmatrix}
-
\begin{bmatrix}
 x_0 & 0 \\
 0   & x_0 
\end{bmatrix}
)
=0, \\
 x_0 = t_1 -(t_1 + t_2) t_1, \qquad x_0 = t_2 -t_1 t_2. 
\end{gather*}
One easily sees that
$
t_1 = \frac{1}{3},\quad 
t_2 = \frac{2}{9}
$ and hence 
$
x_0 = \frac{4}{27}.
$
We will next compute $c_n^{G_3}$ by using Theorem 8.1.
Put 
$
w_i(x) = f_i(x)-1, i=1,2
$
and
$$
F_1(w_1,w_2) = (w_1 +1) \{(w_1 +1) + (w_2 +1)\},
\qquad
F_2(w_1,w_2) = (w_1 +1)  (w_2 +1).
$$
We then have
$$
w_i(x) = x F_i(w_1(x), w_2(x)), \qquad i=1,2.
$$
It follows that by (10.1)
\begin{gather*}
F_2^n(w_1(z),w_2(z))  = (w_2(z) +1 )^{3n}, \\
w'_1(z)  =2 ( w_2(z) + 1) w'_2(z).
\end{gather*}
By Theorem 8.1, one has
\begin{align*}
c_n^{G_3}(1)
& =\frac{1}{n}
\frac{1}{2\pi \sqrt{-1}} \int_C \! \!
\frac{ F_2^n(w_1(z),w_2(z))}{w_2^n(z)} w'_1 (z) dz \\
& =\frac{1}{n}
\frac{1}{2\pi \sqrt{-1}} \int_C \! \!
\frac{(w_2(z) +1 )^{3n}}{w_2^n(z)} 2 ( w_2(z) + 1) w'_2(z) dz \\
& =\frac{2}{n}
\frac{1}{2\pi \sqrt{-1}} \int_C \! \!
\frac{ (w_2 +1)^{3n +1}}{w_2^n} dw_2 \\
& =\frac{2}{n}
\frac{1}{2\pi \sqrt{-1}} 
\sum_{k=0}^{3n+1}
\binom{3n+1}{k}
\int_C \! \!
w_2^{k-n} dw_2 \\
& = \frac{2}{n} \binom{3n+1}{n-1}.
\end{align*}
Similarly we have
\begin{align*}
c_n^{G_3}(2)
& =\frac{1}{n}
\frac{1}{2\pi \sqrt{-1}} \int_C \! \!
\frac{ F_2^n(w_1(z),w_2(z))}{w_2^n(z)} w'_2 (z) dz \\
& =\frac{1}{n}
\frac{1}{2\pi \sqrt{-1}} \int_C \! \!
\frac{(w_2(z) +1 )^{3n}}{w_2^n(z)}  w'_2(z) dz \\
& =\frac{1}{n}
\frac{1}{2\pi \sqrt{-1}} \int_C \! \!
\frac{ (w_2 +1)^{3n }}{w_2^n} dw_2 \\
& =\frac{1}{n}
\frac{1}{2\pi \sqrt{-1}} 
\sum_{k=0}^{3n}
\binom{3n}{k}
\int_C \! \!
w_2^{k-n} dw_2 \\
& = \frac{1}{n} \binom{3n}{n-1}.
\end{align*}
Thereofore we obtain
$$ 
c_n^{G_3}=c_n^{G_3}(1) + c_n^{G_3}(2) = 
\frac{2}{n} \binom{3n+1}{n-1} + \frac{1}{n} \binom{3n}{n-1}
=\frac{2}{n+1} \binom{3n}{n}.
$$
Hence 
$$
c_0^{G_3}=2,\quad
c_1^{G_3}=3,\quad
c_2^{G_3}=10,\quad
c_3^{G_3}=42,\quad
c_4^{G_3}=198,\quad \dots.
$$
This sequence is regarded to be the Fibonacci version of the Catalan numbers.

We remark the following equality on the radius $R_{G_3}$ of convergence of the function $f^{G_3}(x)$:
 \begin{equation*}
 \lim_{n \to \infty} \frac{c_n^G}{c_{n+1}^G}
  =\lim_{n \to \infty} \frac{c_n^G(2)}{c_{n+1}^G(2)} 
  = \lim_{n \to \infty} 
\frac{(n+1) 3 n ! (2n+3) (2n+2)}{(3n+3)!}
 = \frac{4}{27}.
 \end{equation*}

We note that 
by the equality
$f_1(x) = f_2(x)^2$ one has
$$
c_n^{G_3}(1) = \sum_{k=0}^n c_k^{G_3}(2)c_{n-k}^{G_3}(2)
$$
so that the formula 
$$
 \frac{2}{n} \binom{3n+1}{n-1} 
 = \sum_{k=0}^n \frac{1}{k}\binom{3k}{k-1} \frac{1}{n-k}\binom{3n-3k}{n-k-1}
 $$
 holds.

We denote by $\beta$ the golden ratio $\frac{1 + \sqrt{5}}{2}$.
It is the Perron-Frobenius eigenvalue of the matrix
$A^{G_3} =
\begin{bmatrix}
1 & 1 & 0 \\
0 & 0 & 1 \\
1 & 1 & 0 \\ 
\end{bmatrix}.
$
The vector
$\begin{bmatrix}
\beta^{-2}\\
\beta^{-3}\\
\beta^{-2}\\
\end{bmatrix}
$
is the unique positive eigenvector for the eigenvalue $\beta$ whose sum is one.
It then folllows that
$\varphi(P_1) = \beta^{-1}, \varphi(P_2) = \beta^{-2}.
$
Therefore we have
\begin{align*}
c^{G_3,\varphi}_n (1) 
& = c^{G_3}_n (1)\cdot \varphi(P_1)
  = \frac{1}{\beta}\frac{2}{n} \binom{3n+1}{n-1},\\  
c^{G_3,\varphi}_n (1) 
& = c^{G_3}_n (2)\cdot \varphi(P_2)
  = \frac{1}{\beta^2}\frac{1}{n} \binom{3n}{n-1}  
\end{align*}
and
$$
c^{G_3,\varphi}_n =c^{G_3,\varphi}_n (1) + c^{G_3,\varphi}_n (2)
= \frac{2n\beta +(1-n)}{(n+1)(2n+1)}\binom{3n}{n}
= \frac{n\sqrt{5} +1}{(n+1)(2n+1)}\binom{3n}{n}.
$$

\section{Concluding remarks}
Let $\Lambda_G$ be the topological Markov shift of the edges of the graph $G$ 
as in Section 3.
Let $F(A^G)$ be the sub-Fock space associated with the matrix $A^G$.
It is the Hilbert space of the direct sum 
$F(A^G) = \bigoplus_{n=0}^\infty F_n(A^G)$
of the sequence 
$F_n(A^G)$ of the finite dimensional Hilbert spaces
whose orthonomal basis consists of the vectors indexed by the admissible words of the topological Markov shift $\Lambda_G$ of length $n$ for $n =1,2,\dots$.
For $n=0$,
the space $F_0(A^G)$ denotes the one dimensional vector space 
${\Bbb C}\Omega$ of the vacuum vector $\Omega$.
Let $T_e, e \in E$ be the creation operators on $F(A^G)$.
By \cite{EFW81} (cf.\cite{Ev82}),
the quotient images $\bar{T}_e, e\in E$ by the $C^*$-algebra of compact operators on $F(A^G)$ satisfy the relations (3.1).  
We put
$T_G = \sum_{e \in E} T_e$.
 Then the equality
 $<(T_G + T_G^*)^{2n} \Omega \mid \Omega > = c_n^G
 $
 was pointed out by Yoshimichi Ueda.
 The author thanks to him for his suggestion.
Related discussions are seen in several papers of free probability theory 
(cf. \cite{Dyk},\cite{HP},\cite{Vo86},\cite{Vo90}, etc.).   
 
By the above formula of the $G$-Catalan numbers, one may generalize the numbers
$c_n^G$ to general subshifts. 
For a general subshift $\Lambda$ over alphabet $\Sigma$, 
let $F(\Lambda)$ be the sub-Fock Hilbert space associated with it (\cite{Ma}).
  The creation operators $T_\alpha$ for $\alpha \in \Sigma$
  are similarly defined.
  We put
  the operator
$T_\Lambda = \sum_{\alpha \in \Sigma} T_\alpha$ on  $F(\Lambda)$.
We may define the $\Lambda$-Catalan numbers $c^\Lambda_n$ by the formula
$$
c^\Lambda_n = <(T_\Lambda + T_\Lambda^*)^{2n} \Omega \mid \Omega >.
 $$
 The numbers will be studied in \cite{Ma10}

\medskip

{\it Acknowledgments:}\, 
The author would like to thank Toshihiro Hamachi, Wolfgang Krieger and Yoshimichi Ueda for their fruitful suggestions and comments.
Yoshimichi Ueda also kindly informed to the author the Land's text book \cite{La} and some references on this subject. 
The author also thank Hiroaki Yoshida who kindly informed to the author the Deutsch's paper \cite{Deu}.

\end{document}